\documentclass[12pt,a4paper]{wang}%{amsart}
\usepackage{amssymb,amstext,amsmath,amscd,amsfonts,mathdots}
\usepackage{amsmath}% for the equation label
\numberwithin{equation}{section}
\usepackage{etex}% for the counters in tex
\usepackage[numbers]{natbib}% Citations using numbers
\usepackage{hyperref}% for ref
\usepackage{amsthm}% for proof
\usepackage{setspace}% for spacing
\usepackage[vcentermath]{youngtab}%for young diagram
\usepackage{ytableau}
\usepackage[all]{xy}% for quiver
\usepackage{tikz}% for tikz quiver
\usetikzlibrary{shapes.geometric, arrows}
\usetikzlibrary{calc}
\usepackage{enumitem}% for [(1)]
\setlist[description]{leftmargin=\parindent,labelindent=\parindent}
\usepackage{diagbox}% for diag
\usepackage{multirow,multicol}
\usepackage{pdflscape}% for landscape
\usepackage{indentfirst}%for indent
\usepackage{color}% for color
\usepackage{colortbl}
\definecolor{mygray}{gray}{0.80}
\definecolor{myblue}{rgb}{0.50,0.50, 0.95}
\definecolor{myred}{rgb}{0.99, 0.51, 0.65}
\usepackage{lipsum}% for the bib spacing
\let\OLDthebibliography\thebibliography
\renewcommand\thebibliography[1]{\OLDthebibliography{#1}
\setlength{\parskip}{0pt}
\setlength{\itemsep}{0.5pt}}
\newtheorem{theorem}{Theorem}[section]
\newtheorem{theorem*}{Theorem}
\newtheorem{corollary}[theorem]{Corollary}
\newtheorem{corollary*}[theorem*]{Corollary}
\newtheorem{lemma}[theorem]{Lemma}
\newtheorem{proposition}[theorem]{Proposition}

\theoremstyle{definition}
\newtheorem{definition}[theorem]{Definition}
\newtheorem{remark}[theorem]{Remark}

\newtheorem*{question*}{Question}

\newtheorem*{conjecture*}{Conjecture}

\newtheorem*{notation*}{Notation}

\newtheorem*{claim*}{Claim}
\newtheorem{definitiontheorem}[theorem]{Definition-Theorem}
\begin{document}
\begin{spacing}{1.2}

\title{\textbf{$\tau$-tilting finiteness of two-point algebras I}}
\author{Qi Wang}

\abstract{\ \ \ \ As the first attempt to classify $\tau$-tilting finite two-point algebras, we have determined the $\tau$-tilting finiteness for minimal wild two-point algebras and some tame two-point algebras.}

\maketitle

%%%%%%%%%%%%%%%%%%%%%%%%%%%%%%%%%%%%%%%%%%%%%%%%%%%%%%%%%%%%%%%%%%%%%%%%%%%%%%%%%%%%%%%%%%%%%%%%%%%%%%%%%%%%%%%%%%%%%%%%%%%%%%
\section{Introduction}
Throughout this paper, we always assume that $\Lambda$ is a finite-dimensional basic algebra over an algebraically closed field $K$. In particular, the representation type of $\Lambda$ is divided into representation-finite, (infinite-)tame and wild.

$\tau$-tilting theory is introduced by Adachi, Iyama and Reiten \cite{AIR}, in which they constructed support $\tau$-tilting modules as a generalization of the classical tilting modules. We recall that a right $\Lambda$-module $M$ is called support $\tau$-tilting if $\mathsf{Hom}_{\Lambda_M}(M,\tau M)=0$ and $\left | M \right |=\left | \Lambda_M \right |$ taking over $\Lambda_M:=\Lambda/\Lambda (1-e)\Lambda$, namely, $e$ is an idempotent of $\Lambda$ such that the simple summands of $e\Lambda/(e \mathsf{rad}\ \Lambda)$ are exactly the simple composition factors of $M$. Moreover, a support $\tau$-tilting $\Lambda$-module $M$ is called $\tau$-tilting if $\Lambda_M=\Lambda$. This wider class bijectively corresponds to the class of two-term silting complexes, functorially finite torsion classes, left finite semibricks and so on. We refer to \cite{AIR} and \cite{Asai} for details.

We are interested in $\tau$-tilting finite algebras studied in \cite{DIJ-tau-tilting-finite}, that is, algebras with only finitely many pairwise non-isomorphic basic (support) $\tau$-tilting modules. It is obvious that a representation-finite algebra is $\tau$-tilting finite. Also, it is not difficult to find a tame or a wild algebra which is $\tau$-tilting finite. The $\tau$-tilting finiteness for several classes of algebras has been determined, such as algebras with radical square zero \cite{Ada-rad-square-0}, preprojective algebras of Dynkin type \cite{Mizuno}, Brauer graph algebras \cite{AAC-brauer grapha}, biserial algebras \cite{Mousavand-biserial alg} and classical Schur algebras \cite{W-schur}. In particular, it has been proved in some cases that $\tau$-tilting finiteness coincides with representation-finiteness, including gentle algebras \cite{P-gentle}, cycle-finite algebras \cite{MS-cycle-finite}, tilted and cluster-tilted algebras \cite{Z-tilted}, simply connected algebras \cite{W-simply}, quasi-tilted algebras, locally hereditary algebras, etc. \cite{AHMW-joint}.

We notice that local algebras, i.e., algebras with only one simple module, are always $\tau$-tilting finite. This motivates us to study the algebras with exactly two simple modules (up to isomorphism), which are called \emph{two-point algebras} in this paper. We point out that Aihara-Kase \cite{AK-two-point} and Kase \cite{K-two-point} have got some interesting results. For example, Kase \cite[Theorem 6.1]{K-two-point} showed that one can always find a $\tau$-tilting finite two-point algebra such that the Hasse quiver of the poset of pairwise non-isomorphic basic support $\tau$-tilting modules, is isomorphic to a $t$-gon ($t\geqslant 4$). Besides, it is well-known that Kronecker algebra $K(\xymatrix@C=0.7cm{\bullet\ar@<0.5ex>[r]^{\ }\ar@<-0.5ex>[r]_{\ }&\bullet})$ is $\tau$-tilting infinite (we present a proof in Lemma \ref{kronecker} for the convenience of readers). Thus, $\Lambda$ is $\tau$-tilting infinite if the quiver of $\Lambda$ contains multiple arrows.

It is worth mentioning that two-point algebras are fundamental if we consider the representation type of general algebras, and the representation type of two-point algebras has been determined for many years. We may review these results here: the maximal representation-finite two-point algebras are classified by Bongartz and Gabriel \cite{BG-covering}, the tame two-point algebras are classified by several authors in \cite{BH-two-point-without-loop,DG-certain-two-point,Geiss-tame-two-point,Han-wild-two-point}, and the minimal wild two-point algebras are classified by Han \cite{Han-wild-two-point}.

As we mentioned above, $\tau$-tilting finiteness is related to representation type in some classes of algebras. Thus, in order to find a complete classification of $\tau$-tilting finite two-point algebras, it will be useful to determine the $\tau$-tilting finiteness of all minimal wild two-point algebras. We recall that a complete list of minimal wild two-point algebras is given by Han \cite{Han-wild-two-point}, which is displayed by Table W in his paper. (See also Appendix A of this paper.) Then, the first main result in this paper is presented as follows.

\begin{theorem}\label{table-W}
Let $W_i$ be a minimal wild two-point algebra from Table W. Then,
\begin{description}\setlength{\itemsep}{-3pt}
  \item[(1)] $W_1$, $W_2$, $W_3$ and $W_5$ are $\tau$-tilting infinite.
  \item[(2)] $W_4$ and $W_6\sim W_{34}$ are $\tau$-tilting finite. Moreover, we have
  \begin{center}
\renewcommand\arraystretch{1.2}
\begin{tabular}{|c|c|c|c|c|c|c|c|c|c|c|c|}
\hline
$W_i$&$W_4$&$W_6 $&$ W_7 $&$ W_8$&$W_9 $&$ W_{10}$&$ W_{11}$&$W_{12}$&$W_{13}$&$W_{14}$ \\ \hline
$\#\mathsf{s\tau\text{-}tilt}\ W_i$&5&6&\multicolumn{2}{c|}{8} &\multicolumn{2}{c|}{6} &7&\multicolumn{2}{c|}{5}&10\\ \hline
Type&$\mathcal{H}_{1,2}$&$\mathcal{H}_{1,3}$&\multicolumn{2}{c|}{$\mathcal{H}_{1,5}$}&\multicolumn{2}{c|}{$\mathcal{H}_{1,3}$}&$\mathcal{H}_{1,4}$&\multicolumn{2}{c|}{$\mathcal{H}_{1,2}$}&$\mathcal{H}_{3,5}$ \\ \hline\hline
$W_i$&$W_{15} $&$ W_{16}$&$W_{17} $&$W_{18} $&$ W_{19} $&$ W_{20} $&$ W_{21} $&$W_{22} $&$ W_{23} $&$ W_{24}$ \\   \hline
$\#\mathsf{s\tau\text{-}tilt}\ W_i$&9&8&9&8&\multicolumn{4}{c|}{7}&8&10 \\ \hline
Type&$\mathcal{H}_{2,5}$&$\mathcal{H}_{3,3}$&$\mathcal{H}_{2,5}$&$\mathcal{H}_{3,3}$&\multicolumn{4}{c|}{$\mathcal{H}_{2,3}$}&$\mathcal{H}_{3,3}$&$ \mathcal{H}_{3,5}$  \\ \hline\hline
$W_i$&$W_{25} $&$ W_{26}$&$W_{27} $&$W_{28} $&$ W_{29} $&$ W_{30} $&$ W_{31} $&$W_{32}$&$ W_{33} $&$ W_{34}$ \\   \hline
$\#\mathsf{s\tau\text{-}tilt}\ W_i$&\multicolumn{2}{c|}{7}&\multicolumn{4}{c|}{8}&\multicolumn{4}{c|}{6} \\ \hline
Type&\multicolumn{2}{c|}{$\mathcal{H}_{2,3}$}&$\mathcal{H}_{2,4}$&\multicolumn{2}{c|}{$\mathcal{H}_{3,3}$}&$\mathcal{H}_{2,4}$&\multicolumn{4}{c|}{$\mathcal{H}_{2,2}$}  \\ \hline
\end{tabular},
\end{center}
where $\#\mathsf{s\tau\text{-}tilt}\ W_i$ is the number of isomorphism classes of basic support $\tau$-tilting $W_i$-modules and the type of $\mathcal{H}(\mathsf{s\tau\text{-}tilt}\ W_i)$ is defined in Definition \ref{type}.
\end{description}
\end{theorem}

According to Theorem \ref{table-W}, most of the minimal wild two-point algebras are $\tau$-tilting finite so that one may expect to give a complete classification of $\tau$-tilting finite tame two-point algebras. This should also be useful toward the complete classification of $\tau$-tilting finite two-point algebras. However, it is difficult at this moment to give a complete result on tame two-point algebras, because the tameness of two-point algebras depends on the technique called \emph{degeneration}, and it is still open to finding the relation between $\tau$-tilting finiteness and degeneration.

We may give a partial result on tame two-point algebras. We recall from \cite{Han-wild-two-point} (see also Proposition \ref{han-two-point} of this paper) that all tame two-point algebras can degenerate to a finite set (Table T in \cite{Han-wild-two-point}) of two-point algebras. Then, we check the $\tau$-tilting finiteness for algebras in Table T\footnote{We mention that some relations are omitted in the original Table T in \cite{Han-wild-two-point} so that several algebras (e.g., $T_4$ and $T_5$) in the original Table T are not finite-dimensional. However, we have added these omitted relations in this paper so that all algebras in Table T are finite-dimensional, see Appendix A.} as follows. This is the second main result in this paper.

\begin{theorem}\label{table-T}
Let $T_i$ be an algebra from Table T.
\begin{description}\setlength{\itemsep}{-3pt}
  \item[(1)] $T_1$, $T_3$ and $T_{17}$ are $\tau$-tilting infinite.
  \item[(2)] Others are $\tau$-tilting finite. Moreover, we have the following posets,
\begin{center}
\renewcommand\arraystretch{1.2}
\begin{tabular}{|c|c|c|c|c|c|c|c|c|c|}
\hline
$T_i$&$T_2$&$T_4$&$T_5$&$T_6$&$T_7$&$T_8$&$T_9$&$T_{10}$&$T_{11}$  \\ \hline
$\#\mathsf{s\tau\text{-}tilt}\ T_i$&\multicolumn{2}{c|}{6} &5&6&\multicolumn{2}{c|}{5} &8&12&8 \\ \hline
Type&\multicolumn{2}{c|}{$\mathcal{H}_{1,3}$} &$\mathcal{H}_{1,2}$&$\mathcal{H}_{1,3}$&\multicolumn{2}{c|}{$\mathcal{H}_{1,2}$} &$\mathcal{H}_{3,3}$&$\mathcal{H}_{5,5}$&$\mathcal{H}_{3,3}$ \\ \hline\hline
$T_i$&$T_{12}$ &$ T_{13} $ &$ T_{14} $ &$ T_{15}$ &$ T_{16}$ &$T_{18}$ &$ T_{19}$ &$ T_{20} $ &$ T_{21}$ \\   \hline
$\#\mathsf{s\tau\text{-}tilt}\ T_i$&7&6&8&7&9 &8&6&7&6  \\ \hline
Type&$\mathcal{H}_{2,3}$&$\mathcal{H}_{2,2}$&$\mathcal{H}_{3,3}$&$\mathcal{H}_{2,3}$ &$\mathcal{H}_{2,5}$ &$\mathcal{H}_{3,3}$&$\mathcal{H}_{2,2}$&$\mathcal{H}_{2,3}$&$\mathcal{H}_{2,2}$  \\ \hline
\end{tabular}.
\end{center}
\end{description}
\end{theorem}

We observe that Theorem \ref{table-W} and Theorem \ref{table-T} are useful to determine the $\tau$-tilting finiteness for several other classes of algebras, such as tame two-point distributive algebras \cite{Geiss-tame-two-point}, two-point symmetric special biserial algebras \cite{AIP-joint} and so on. We have given an easy observation in Proposition \ref{application}.

This paper is organized as follows. In Section 2, we review some basic concepts of $\tau$-tilting theory and silting theory. Besides, we list some reduction theorems that we will use and carry out several explicit computations. In Section 3, we give the proofs of Theorem \ref{table-W} and Theorem \ref{table-T}.

\vspace{0.5cm}

\noindent\textit{Acknowledgements.} I am grateful to my supervisor Prof. Susumu Ariki for suggesting this project and for his kind guidance during the writing of this paper. A lot of thanks to Dr. Kengo Miyamoto for teaching me numerous knowledge of $\tau$-tilting theory, to Prof. Takuma Aihara and Prof. Ryoichi Kase for many useful discussions. I appreciate the anonymous referee very much for being patient with me and giving me many useful suggestions. Finally, I would like to mention my sincere gratitude to Prof. Toshiaki Shoji. Without his recommendation, I could not imagine studying abroad.

\section{Preliminaries}
We refer to \cite{ASS} for the background materials on the representation theory of finite-dimensional algebras and the basic knowledge of quiver representations.

Let $\mathsf{mod}\ \Lambda$ be the category of finitely generated right $\Lambda$-modules and $\mathsf{proj}\ \Lambda$ the full subcategory of $\mathsf{mod}\ \Lambda$ consisting of projective $\Lambda$-modules. For any $M\in \mathsf{mod}\ \Lambda$, we denote by $\mathsf{add}(M)$ (respectively, $\mathsf{Fac}(M)$) the full subcategory of $\mathsf{mod}\ \Lambda$ whose objects are direct summands (respectively, factor modules) of finite direct sums of copies of $M$. We often describe $\Lambda$-modules via their composition series. For example, each simple module $S_i$ is written as $i$, and then $\substack{1\\2}=\substack{S_1\\S_2}$ is an indecomposable $\Lambda$-module $M$ with a unique simple submodule $S_2$ such that $M/S_2\simeq S_1$.

We denote by $\Lambda^{\mathsf{op}}$ the opposite algebra of $\Lambda$ and by $|M|$ the number of isomorphism classes of indecomposable direct summands of $M$. Let $\tau$ be the Auslander-Reiten translation on the module category. Note that it is not functorial.
\begin{definition}(\cite[Definition 0.1]{AIR})
Let $M\in \mathsf{mod}\ \Lambda$. Then,
\begin{description}\setlength{\itemsep}{-3pt}
  \item[(1)] $M$ is called $\tau$-rigid if $\mathsf{Hom}_\Lambda(M,\tau M)=0$.
  \item[(2)] $M$ is called $\tau$-tilting if $M$ is $\tau$-rigid and $\left | M \right |=\left | \Lambda \right |$.
  \item[(3)] $M$ is called support $\tau$-tilting if there exists an idempotent $e$ of $\Lambda$ such that $M$ is a $\tau$-tilting $\left ( \Lambda/\Lambda e \Lambda\right )$-module.
\end{description}
\end{definition}

Let $(M,P)$ be a pair with $M\in \mathsf{mod}\ \Lambda$ and $P\in \mathsf{proj}\ \Lambda$. We recall that $(M,P)$ is said to be a support $\tau$-tilting pair if $M$ is $\tau$-rigid, $\mathsf{Hom}_\Lambda(P,M)=0$ and $|M|+|P|=|\Lambda|$. This is actually an equivalent definition for support $\tau$-tilting modules, i.e., $(M,P)$ is a support $\tau$-tilting pair if and only if $M$ is a $\tau$-tilting $\left ( \Lambda/\Lambda e\Lambda \right )$-module and $P=e\Lambda$.

We denote by $\mathsf{\tau\text{-}rigid}\ \Lambda$ (respectively, $\mathsf{s\tau\text{-}tilt}\ \Lambda$) the set of isomorphism classes of indecomposable $\tau$-rigid (respectively, basic support $\tau$-tilting) $\Lambda$-modules. It is known from \cite{AIR} that any $\tau$-rigid $\Lambda$-module is a direct summand of some $\tau$-tilting $\Lambda$-module.
\begin{definition}
An algebra $\Lambda$ is called $\tau$-tilting finite if it has only finitely many pairwise non-isomorphic basic $\tau$-tilting modules. Otherwise, $\Lambda$ is called $\tau$-tilting infinite.
\end{definition}

Let $\mathcal{C}$ be an additive category and $X, Y$ objects of $\mathcal{C}$. A morphism $f: X\rightarrow Z$ with $Z\in \mathsf{add}(Y)$ is called a minimal left $\mathsf{add}(Y)$-approximation of $X$ if it satisfies:
\begin{itemize}\setlength{\itemsep}{-3pt}
  \item every $h\in \mathsf{Hom}_\mathcal{C}(Z,Z)$ that satisfies $h\circ f=f$ is an automorphism,
  \item $\mathsf{Hom}_\mathcal{C}(f,Z'): \mathsf{Hom}_\mathcal{C}(Z,Z')\longrightarrow \mathsf{Hom}_\mathcal{C}(X,Z')$ is surjective for any $Z'\in \mathsf{add}(Y)$,
\end{itemize}
where $\mathsf{add}(Y)$ is the category of all direct summands of finite direct sums of copies of $Y$.

We recall the concept of left mutation which is the core of $\tau$-tilting theory.
\begin{definitiontheorem}(\cite[Definition 2.19, Theorem 2.30]{AIR})
Let $T=M\oplus N$ be a basic support $\tau$-tilting $\Lambda$-module with an indecomposable direct summand $M$ satisfying $M\notin \mathsf{Fac}(N)$. We take an exact sequence with a minimal left $\mathsf{add}(N)$-approximation $f$:
\begin{center}
$M\overset{f}{\longrightarrow}N'\longrightarrow \mathsf{coker}\ f \longrightarrow 0$.
\end{center}
We call $\mu_M^-(T):=(\mathsf{coker}\ f)\oplus N$ the left mutation of $T$ with respect to $M$, which is again a basic support $\tau$-tilting $\Lambda$-module. (The right mutation $\mu_M^+(T)$ can be defined dually.)
\end{definitiontheorem}

In the above, Zhang has shown in \cite{Z-mutation} that $\mathsf{coker}\ f$ is either $0$ or indecomposable. Moreover, one can show that $\mathsf{coker}\ f$ cannot be projective.

We may construct a directed graph $\mathcal{H}(\mathsf{s\tau\text{-}tilt}\ \Lambda)$ by drawing an arrow from $T_1$ to $T_2$ if $T_2$ is a left mutation of $T_1$. On the other hand, we can regard $\mathsf{s\tau\text{-}tilt}\ \Lambda$ as a poset with respect to a partial order $\leq $. For any $M, N\in \mathsf{s\tau\text{-}tilt}\ \Lambda$, we say that $N\leq M$ if $\mathsf{Fac}(N) \subseteq \mathsf{Fac}(M)$. Then, the directed graph $\mathcal{H}(\mathsf{s\tau\text{-}tilt}\ \Lambda)$ is exactly the Hasse quiver on the poset $\mathsf{s\tau\text{-}tilt}\ \Lambda$, see \cite[Corollary 2.34]{AIR}.

The following statement implies that an algebra $\Lambda$ is $\tau$-tilting finite if we can find a finite connected component of $\mathcal{H}(\mathsf{s\tau\text{-}tilt}\ \Lambda)$.
\begin{proposition}{\rm{(\cite[Corollary 2.38]{AIR})}} \label{finite-connected-component}
If the Hasse quiver $\mathcal{H}(\mathsf{s\tau\text{-}tilt}\ \Lambda)$ contains a finite connected component $\Delta$, then $\mathcal{H}(\mathsf{s\tau\text{-}tilt}\ \Lambda)=\Delta$.
\end{proposition}

\subsection{Silting theory}
We denote by $\mathsf{C^b(proj}\ \Lambda)$ the category of bounded complexes of projective $\Lambda$-modules and by $\mathsf{K^b(proj}\ \Lambda)$ the corresponding homotopy category which is triangulated. Besides, we denote by $\sim_h$ the homotopy equivalence in $\mathsf{K^b(proj}\ \Lambda)$. For any $T\in \mathsf{K^b(proj}\ \Lambda)$, let $\mathsf{thick}\ T$ be the smallest full subcategory of $\mathsf{K^b(proj}\ \Lambda)$ containing $T$, which is closed under cones, $[\pm1]$, direct summands and isomorphisms.
\begin{definition}(\cite[Definition 2.1]{AI-silting})
A complex $T\in \mathsf{K^b(proj}\ \Lambda)$ is called presilting if
\begin{center}
$\mathsf{Hom}_{\mathsf{K^b(proj}\ \Lambda)}(T,T[i])=0$ for any $i>0$.
\end{center}
A presilting complex $T$ is called silting if $\mathsf{thick}\ T=\mathsf{K^b(proj}\ \Lambda)$.
\end{definition}

Similar to the left mutation of support $\tau$-tilting modules, we recall the irreducible left silting mutation (\cite[Definition 2.30]{AI-silting}) of silting complexes. Let $T=X\oplus Y$ be a basic silting complex in $\mathsf{K^b(proj}\ \Lambda)$ with an indecomposable summand $X$. We take a minimal left $\mathsf{add}(Y)$-approximation $\pi$ and a triangle
\begin{center}
$X\overset{\pi}{\longrightarrow}Z\longrightarrow \mathsf{cone}(\pi)\longrightarrow X[1]$.
\end{center}
Then, $\mathsf{cone}(\pi)$ is indecomposable and $\mu_X^-(T):= \mathsf{cone}(\pi)\oplus Y$ is again a basic silting complex in $\mathsf{K^b(proj}\ \Lambda)$, see \cite[Theorem 2.31]{AI-silting}. We call $\mu_X^-(T)$ the irreducible left mutation of $T$ with respect to $X$.

A complex in $\mathsf{K^b(proj}\ \Lambda)$ is called \emph{two-term} if it is homotopy equivalent to a complex $T$ concentrated in degree $0$ and $-1$, i.e.,
\begin{center}
$( T^{-1}\overset{d_T^{-1}}{\longrightarrow } T^0 ):=\xymatrix@C=0.7cm{\ldots\ar[r]&0\ar[r]&T^{-1}\ar[r]^{d_T^{-1}}&T^0\ar[r]&0\ar[r]&\ldots}$.
\end{center}
We denote by $\mathsf{2\text{-}silt}\ \Lambda$ the set of isomorphism classes of basic two-term silting complexes in $\mathsf{K^b(proj}\ \Lambda)$. Similarly, there is a partial order $\leq$ on the set $\mathsf{2\text{-}silt}\ \Lambda$ which is introduced by \cite[Theorem 2.11]{AI-silting}. For any $S, T\in\mathsf{2\text{-}silt}\ \Lambda$, we say that $S\leq T$ if $\mathsf{Hom}_{\mathsf{K^b(proj}\ \Lambda)}(T,S[i])=0$ for any $i>0$. We denote by $\mathcal{H}(\mathsf{2\text{-}silt}\ \Lambda)$ the Hasse quiver of $\mathsf{2\text{-}silt}\ \Lambda$, which is compatible with the irreducible left mutations of silting complexes.

\begin{proposition}{\rm{(\cite[Theorem 3.2]{AIR})}}\label{tilting-silting}
There exists a poset isomorphism between $\mathsf{s\tau\text{-}tilt}\ \Lambda$ and $2\mathsf{\text{-}silt}\ \Lambda$. More precisely, the bijection is given by
\begin{center}
$\xymatrix@C=0.8cm@R=0.2cm{M\ar@{|->}[r]&(P_1\oplus P\overset{[f,0]}{\longrightarrow} P_0)}$
\end{center}
where $(M,P)$ is the corresponding support $\tau$-tilting pair and $P_1\overset{f}{\longrightarrow }P_0\overset{}{\longrightarrow }M\longrightarrow 0$ is a minimal projective presentation of $M$.
\end{proposition}

Since the poset $\mathsf{s\tau\text{-}tilt}\ \Lambda$ has the unique maximal element $\Lambda$ and the unique minimal element $0$, we can define the type of $\mathcal{H}(\mathsf{s\tau\text{-}tilt}\ \Lambda)$ (equivalently, $\mathcal{H}(\mathsf{2\mathsf{\text{-}silt}}\ \Lambda)$) as follows.
\begin{definition}\label{type}
Let $\Lambda$ be a $\tau$-tilting finite algebra. We say that the Hasse quiver $\mathcal{H}(\mathsf{s\tau\text{-}tilt}\ \Lambda)$ is of type $\mathcal{H}_{m,n}$ if it is of the form
\begin{center}
$\xymatrix@C=0.7cm@R=0.02cm{&\triangle_1 \ar[r]&\triangle_2\ar[r] &\ldots \ar[r]&\triangle_m \ar[dr]&\\
\Lambda \ar[ur]\ar[dr]&&&&&0.\\
&\square_1 \ar[r]&\square_2 \ar[r]&\ldots \ar[r]&\square_n \ar[ur]&} $
\end{center}
Moreover, we have $\mathcal{H}_{m,n}\simeq \mathcal{H}_{n,m}$.
\end{definition}

We have the following equivalent condition for $\Lambda$ to be $\tau$-tilting finite.
\begin{proposition}{\rm{(\cite[Corollary 2.9]{DIJ-tau-tilting-finite})}}\label{tau-tilting-finite-rigid}
An algebra $\Lambda$ is $\tau$-tilting finite if and only if one of (equivalently, any of) the sets $\mathsf{\tau\text{-}rigid}\ \Lambda$, $\mathsf{s\tau\text{-}tilt}\ \Lambda$ and $\mathsf{2\text{-}silt}\ \Lambda$ is finite.
\end{proposition}

\subsection{Reduction theorems}
There are some reduction theorems. First, we review the brick-$\tau$-rigid correspondence introduced by Demonet, Iyama and Jasso \cite{DIJ-tau-tilting-finite}. Recall that $M\in \mathsf{mod}\ \Lambda$ is called a \textit{brick} if $\mathsf{End}_\Lambda(M)=K$. We denote by $\mathsf{brick}\ \Lambda$ the set of isomorphism classes of bricks in $\mathsf{mod}\ \Lambda$.
\begin{lemma}{\rm{(\cite[Theorem 4.2]{DIJ-tau-tilting-finite})}} \label{rigid-brick}
Let $\Lambda$ be a finite-dimensional algebra. Then, $\Lambda$ is $\tau$-tilting finite if and only if the set $\mathsf{brick}\ \Lambda$ is finite.
\end{lemma}

Let $\Lambda_1$, $\Lambda_2$ be two algebras. We call $\Lambda_2$ a \textit{quotient} or \textit{quotient algebra} of $\Lambda_1$ if there exists a surjective $K$-algebra homomorphism $\phi :\Lambda_1\rightarrow \Lambda_2$.
\begin{corollary}\label{factor}
Suppose that $\Lambda_2$ is a quotient algebra of $\Lambda_1$. If $\Lambda_2$ is $\tau$-tilting infinite, then $\Lambda_1$ is also $\tau$-tilting infinite.
\end{corollary}
\begin{proof}
There exists a $K$-linear fully-faithful functor $\mathcal{F}:\mathsf{mod}\ \Lambda_2\rightarrow \mathsf{mod}\ \Lambda_1$, and $\mathcal{F}$ sends a brick in $\mathsf{mod}\ \Lambda_2$ to a brick in $\mathsf{mod}\ \Lambda_1$. Then, the statement follows from Lemma \ref{rigid-brick}.
\end{proof}

\begin{lemma}{\rm{(\cite[Theorem 2.14]{AIR})}} \label{oppsite algebra}
There exists a poset anti-isomorphism between $\mathsf{s\tau\text{-}tilt}\ \Lambda$ and $\mathsf{s\tau\text{-}tilt}\ \Lambda^{\mathsf{op}}$.
\end{lemma}

\begin{lemma}{\rm{(\cite[Theorem 1]{EJR-reduction})}} \label{centers}
Let $I$ be a two-sided ideal generated by central elements which are contained in the Jacobson radical of $\Lambda$. Then, there exists a poset isomorphism between $\mathsf{s\tau\text{-}tilt}\ \Lambda$ and $\mathsf{s\tau\text{-}tilt}\ (\Lambda/I)$.
\end{lemma}

\begin{lemma}\label{reduced-complex}
If $Y\neq 0$ and
\begin{center}
$T_1:=(\xymatrix@C=1.2cm{0 \ar[r]&X\ar[r]^{\binom{1}{f}\ \ \ \ }&X\oplus Y \ar[r]^{\ \ \ \ (-g\circ f, g)} &Z\ar[r]&0})\in \mathsf{K^b(proj}\ \Lambda)$,

$T_2:=(\xymatrix@C=1.2cm{0 \ar[r]&X\oplus Y \ar[r]^{\left ( \begin{smallmatrix}
f_1 & f_2\\
1 & g\\
h_1 &h_2
\end{smallmatrix} \right )\ \ \ \ } &Z\oplus X\oplus M \ar[r]&0})\in \mathsf{K^b(proj}\ \Lambda)$,
\end{center}
then $T_1\sim_h T_1^r$ and $T_2\sim_h T_2^r$, where
\begin{center}
$T_1^r: =(\xymatrix@C=1cm{0 \ar[r]&Y \ar[r]^{g} &Z\ar[r]&0} )\in \mathsf{K^b(proj}\ \Lambda)$,

$T_2^r: =(\xymatrix@C=1.5cm{0 \ar[r]& Y \ar[r]^{\left ( \begin{smallmatrix}
f_2-f_1\circ g\\
h_2-h_1\circ g
\end{smallmatrix} \right )\ \ \ \ } &Z\oplus M \ar[r]&0})\in \mathsf{K^b(proj}\ \Lambda)$.
\end{center}
\end{lemma}
\begin{proof}
(1) We define $\varphi :T_1\rightarrow T_1^r $ and $\psi :T_1^r\rightarrow T_1$ as follows.
\begin{center}
$\xymatrix@C=1.2cm@R=0.8cm{T_1:&0 \ar[r]&X\ar@<0.5ex>[d]^{0}\ar[r]^{\binom{1}{f}\ \ \ \ }&X\oplus Y \ar@<0.5ex>[d]^{(-f,1)}\ar[r]^{\ \ \ \ (-g\circ f, g)} &Z\ar@<0.5ex>[d]^{1}\ar[r]&0\\
T_1^r:&0\ar[r]&0\ar@{-->}@<0.5ex>[u]^{0}\ar[r]_{0}&Y\ar@{-->}@<0.5ex>[u]^{\binom{0}{1}}\ar[r]_{g} &Z\ar@{-->}@<0.5ex>[u]^{1}\ar[r]&0}$
\end{center}
Then, $\varphi \circ \psi=\text{Id}_{T_1^r}$ and $\psi\circ  \varphi =\left ( 0,\left ( \begin{smallmatrix}
0&0\\
-f&1
\end{smallmatrix} \right ),1 \right )\sim_h \text{Id}_{T_1} $ because the difference $\text{Id}_{T_1}-\psi\circ  \varphi$ is null-homotopic as follows.
\begin{center}
$\xymatrix@C=1.2cm@R=1.2cm{T_1:&0 \ar[r]&X\ar[d]_{1}\ar@{.>}[dl]_{0}\ar[r]^{\binom{1}{f}\ \ \ \ }&X\oplus Y \ar[d]_{\left ( \begin{smallmatrix}
1&0\\
f&0
\end{smallmatrix} \right )}\ar[r]^{\ \ \ \ (-g\circ f, g)}\ar@{.>}[dl]_{(1,0)} &Z\ar@{.>}[dl]_{\ \ \ \ \ \ \ \ 0}\ar[d]_{0}\ar[r]&0\ar@{.>}[dl]_{0}\\
T_1:&0 \ar[r]&X\ar[r]_{\binom{1}{f}\ \ \ \ }&X\oplus Y \ar[r]_{\ \ \ \ (-g\circ f, g)} &Z\ar[r]&0}$
\end{center}

(2) We define $\varphi :T_2\rightarrow T_2^r $ and $\psi :T_2^r\rightarrow T_2$ as follows.
\begin{center}
$\xymatrix@C=1.2cm@R=1cm{T_2:&0 \ar[r]&X\oplus Y \ar@<0.5ex>[d]^{(0,1)}\ar[r]^{\left ( \begin{smallmatrix}
f_1 & f_2\\
1 & g\\
h_1 &h_2
\end{smallmatrix} \right )\ \ \ \ } &Z\oplus X\oplus M \ar@<0.5ex>[d]^{\left ( \begin{smallmatrix}
1 & -f_1&0\\
0 & -h_1&1
\end{smallmatrix} \right )}\ar[r]&0\\
T_2^r:&0 \ar[r]& Y \ar@{-->}@<0.5ex>[u]^{\binom{-g}{1}}\ar[r]_{\left ( \begin{smallmatrix}
f_2-f_1\circ g\\
h_2-h_1\circ g
\end{smallmatrix} \right )\ \ \ \ } &Z\oplus M \ar@{-->}@<0.5ex>[u]^{\left ( \begin{smallmatrix}
1 & 0\\
0 & 0\\
0 & 1
\end{smallmatrix} \right )}\ar[r]&0}$
\end{center}
Then, $\varphi \circ \psi=\text{Id}_{T_2^r}$ and $\psi\circ \varphi =\left ( \left ( \begin{smallmatrix}
0&-g\\
0&1
\end{smallmatrix} \right ),\left ( \begin{smallmatrix}
1&-f_1&0\\
0&0&0\\
0&-h_1&1
\end{smallmatrix} \right )\right )\sim_h \text{Id}_{T_2}$. In fact, the difference $\text{Id}_{T_2}-\psi\circ  \varphi$ is null-homotopic as follows.
\begin{center}
$\xymatrix@C=1.8cm@R=1.3cm{T_2:&0 \ar[r]&X\oplus Y \ar@{.>}[dl]_{0}\ar[d]_{\left ( \begin{smallmatrix}
1 & g\\
0 & 0
\end{smallmatrix} \right )}\ar[r]^{\left ( \begin{smallmatrix}
f_1 & f_2\\
1 & g\\
h_1 &h_2
\end{smallmatrix} \right )\ \ \ \ \  } &Z\oplus X\oplus M \ar@{.>}[dl]_{\left ( \begin{smallmatrix}
0 & 1&0\\
0 & 0&0
\end{smallmatrix} \right )}\ar[d]_{\left ( \begin{smallmatrix}
0 & f_1&0\\
0 & 1&0\\
0&h_1&0
\end{smallmatrix} \right )}\ar[r]&0\ar@{.>}[dl]_{0}\\
T_2:&0 \ar[r]&X\oplus Y \ar[r]_{\left ( \begin{smallmatrix}
f_1 & f_2\\
1 & g\\
h_1 &h_2
\end{smallmatrix} \right )\ \ \ \ \  } &Z\oplus X\oplus M \ar[r]&0}$
\end{center}
Therefore, we have $T_1\sim_h T_1^r$ and $T_2\sim_h T_2^r$.
\end{proof}

\section{Main Results}
In this section, we will prove our main results mentioned in the introduction. But before that, let us review the complete classification for the representation type of two-point algebras. We refer to Appendix A for Table W and Table T.
\begin{proposition}{\rm{(\cite[Theorem 1, Theorem 2]{Han-wild-two-point})}}\label{han-two-point}
Let $\Lambda$ be a two-point algebra. Up to isomorphism and duality, $\Lambda$ is representation-finite or tame if and only if $\Lambda$ degenerates to a quotient algebra of an algebra from Table T, and $\Lambda$ is wild if and only if $\Lambda$ has a minimal wild algebra from Table W as a quotient algebra.
\end{proposition}

We denote by $\mathsf{rad}(\Lambda)$ the Jacobson radical of $\Lambda$ and by $C(\Lambda)$ the center of $\Lambda$. As explained in \cite{EJR-reduction}, although $\Lambda$ has a complicated structure, its quotient algebra
\begin{center}
$\widetilde{\Lambda}:=\Lambda/<C(\Lambda)\cap \mathsf{rad}(\Lambda)>$
\end{center}
may have a simpler structure. Moreover, by Lemma \ref{centers}, we know that $\#\mathsf{s\tau\text{-}tilt}\ \Lambda=\#\mathsf{s\tau\text{-}tilt}\ \widetilde{\Lambda}$ and the Hasse quivers $\mathcal{H}(\mathsf{s\tau\text{-}tilt}\ \Lambda)$ and $\mathcal{H}(\mathsf{s\tau\text{-}tilt}\ \widetilde{\Lambda})$ are of the same type. Then, by using this strategy and Corollary \ref{factor}, we can restrict the algebras (except for $W_4$, $T_{20}$ and $T_{21}$) in Table W and Table T to a small list (i.e., Table $\Lambda$) of two-point algebras. (In Table $\Lambda$, by an algebra $\Lambda_i$, we mean the bound quiver algebra $KQ/I_i$, where $I_i$ is the admissible ideal generated by the relation $(i)$.)

\begin{table}[t]
\begin{center}
\textbf{Table $\Lambda$}\ \ \ \ \ \ \ \ \ \ \ \ \ \
\end{center}
\vspace{-0.3cm}
\rule[0.5\baselineskip]{15cm}{1pt}
\vspace{-0.5cm}
\begin{multicols}{2}
$\Lambda_1=K\left ( \xymatrix@C=1.2cm{1\ar[r]&2}  \right )$;

\rule[0.5\baselineskip]{6cm}{1pt}

$\Lambda_2=K\left ( \xymatrix@C=1.2cm{1\ar@<0.5ex>[r]^{\ }\ar@<-0.5ex>[r]_{\ }&2}  \right )$;

\rule[0.5\baselineskip]{6cm}{1pt}

\ \ \ \ \ $Q: \xymatrix@C=1.2cm{1 \ar[r]^{\mu}& 2\ar@(ur,dr)^{\beta } }$
\vspace{0.3cm}

(3) $\beta^2=0$;

(4) $\beta^3=0$;

\rule[0.5\baselineskip]{6cm}{1pt}

\ \ \ \ \ $Q: \xymatrix@C=1.2cm{1 \ar[r]^{\mu}\ar@(dl,ul)^{\alpha}& 2 \ar@(ur,dr)^{\beta}}$
\vspace{0.3cm}

(5) $\alpha^2=\beta^2=0$;

\rule[0.5\baselineskip]{6cm}{1pt}\\

\ \ \ \ \ $Q: \xymatrix@C=1.2cm{1\ar@<0.5ex>[r]^{\mu}&2\ar@<0.5ex>[l]^{\nu }}$
\vspace{0.3cm}

(6) $\mu\nu=\nu\mu=0$;

\rule[0.5\baselineskip]{6.8cm}{1pt}

\ \ \ \ \ $Q: \xymatrix@C=1.2cm{1 \ar@<0.5ex>[r]^{\mu}\ar@(dl,ul)^{\alpha}& 2 \ar@<0.5ex>[l]^{\nu}}$
\vspace{0.3cm}

(7) $\alpha^2=\mu\nu=\nu\mu=\nu\alpha=0$;

(8) $\alpha^2=\mu\nu=\nu\mu=\nu\alpha\mu=0$;

(9) $\alpha^3=\mu\nu=\nu\mu=\nu\alpha=0$;

(10) $\alpha^3=\mu\nu=\nu\mu=\nu\alpha\mu=\nu\alpha^2\mu=0$;

\rule[0.5\baselineskip]{6.8cm}{1pt}
\end{multicols}

\ \ \ \ \ $Q: \xymatrix@C=1.2cm{1 \ar@<0.5ex>[r]^{\mu}\ar@(dl,ul)^{\alpha}&2 \ar@<0.5ex>[l]^{\nu}\ar@(ur,dr)^{\beta}}$
\vspace{0.3cm}

(11) $\alpha^2=\beta^2=\mu\nu=\nu\mu=\alpha\mu=\beta\nu=0$;

(12) $\alpha^2=\beta^2=\mu\nu=\nu\mu=\beta\nu=\nu\alpha=\alpha\mu\beta=0$.

\rule[0.5\baselineskip]{15cm}{1pt}
\vspace{-0.3cm}
\end{table}

As a preparation for proving Theorem \ref{table-W} and Theorem \ref{table-T}, we shall determine the $\tau$-tilting finiteness of $\Lambda_i$ in Table $\Lambda$. We remark that $\mathcal{H}(\mathsf{s\tau\text{-}tilt}\ \Lambda_1)$ is of type $\mathcal{H}_{1,2}$ and $\mathcal{H}(\mathsf{s\tau\text{-}tilt}\ \Lambda_6)$ is of type $\mathcal{H}_{2,2}$, and we omit the details to show this.

It is well-known that the Kronecker algebra $\Lambda_2$ admits infinitely many (classical) tilting modules, so that it is obviously $\tau$-tilting infinite. Since the construction of tilting modules is not mentioned in this paper, we give a different proof here for the convenience of readers. Similar to the notion of minimal representation-infinite algebra, we call an algebra $\Lambda$ \textit{minimal $\tau$-tilting infinite} if $\Lambda$ is $\tau$-tilting infinite, but any proper quotient algebra of $\Lambda$ is $\tau$-tilting finite.
\begin{lemma}\label{kronecker}
The Kronecker algebra $\Lambda_2$ is minimal $\tau$-tilting infinite.
\end{lemma}
\begin{proof}
It is easy to check that $M_k=\xymatrix@C=0.7cm{K\ar@<0.5ex>[r]^{k }\ar@<-0.5ex>[r]_{1}&K}$ with $k\in K$ is a brick in $\mathsf{mod}\ \Lambda_2$. Since the family $(M_k)_{k\in K}$ consists of infinitely many pairwise non-isomorphic bricks, $\Lambda_2$ is $\tau$-tilting infinite by Lemma \ref{rigid-brick}. Besides, the minimality is obvious.
\end{proof}

\begin{lemma}\label{lambda4}
The two-point algebras $\Lambda_3$ and $\Lambda_4$ are $\tau$-tilting finite.
\end{lemma}
\begin{proof}
Since $\Lambda_3$ is a quotient algebra of $\Lambda_4$, by Corollary \ref{factor}, it suffices to show that $\Lambda_4$ is $\tau$-tilting finite. We show that the poset $\mathsf{2\text{-}silt}\ \Lambda_4$ has a finite connected component and hence, it exhausts all two-term silting complexes in $\mathsf{K^b(proj}\ \Lambda_4)$ by Proposition \ref{finite-connected-component} and Proposition \ref{tilting-silting}. Then, $\Lambda_4$ is $\tau$-tilting finite following from Proposition \ref{tau-tilting-finite-rigid}. Let $P_1$ and $P_2$ be the indecomposable projective $\Lambda_4$-modules. We have
\begin{center}
$P_1=\begin{smallmatrix}
e_1\\
\mu\\
\mu\beta\\
\mu\beta^2
\end{smallmatrix}\simeq \begin{smallmatrix}
1\\
2\\
2\\
2
\end{smallmatrix}$ and $P_2=\begin{smallmatrix}
e_2\\
\beta\\
\beta^2
\end{smallmatrix}\simeq\begin{smallmatrix}
2\\
2\\
2
\end{smallmatrix}$.
\end{center}

We show that $\mathcal{H}$($\mathsf{2\text{-}silt}\ \Lambda_4$) is of type $\mathcal{H}_{1,5}$ as follows,
\begin{center}
\begin{tikzpicture}[shorten >=1pt, auto, node distance=0cm,
   node_style/.style={font=},
   edge_style/.style={draw=black}]
\node[node_style] (v0) at (0,0) {$\left [\begin{smallmatrix}
0\longrightarrow P_1\\
\oplus \\
0\longrightarrow P_2
\end{smallmatrix}  \right ]$};
\node[node_style] (v1) at (12.7,0) {$\left [\begin{smallmatrix}
P_1\longrightarrow 0\\
\oplus \\
0\longrightarrow P_2
\end{smallmatrix}  \right ]$};
\node[node_style] (v) at (12.7,-2) {$\left [\begin{smallmatrix}
P_1\longrightarrow 0\\
\oplus \\
P_2\longrightarrow 0
\end{smallmatrix}  \right ]$};
\node[node_style] (v2) at (0,-2) {$\left [\begin{smallmatrix}
0\longrightarrow P_1\\
\oplus \\
P_2\overset{f_1}{\longrightarrow} P_1^{\oplus 3}
\end{smallmatrix}  \right ]$};
\node[node_style] (v21) at (2.5,-2) {$\left [\begin{smallmatrix}
P_2\overset{f_2}{\longrightarrow} P_1^{\oplus 2}\\
\oplus \\
P_2\overset{f_1}{\longrightarrow} P_1^{\oplus 3}
\end{smallmatrix}  \right ]$};
\node[node_style] (v212) at (5.2,-2) {$\left [\begin{smallmatrix}
P_2\overset{f_2}{\longrightarrow} P_1^{\oplus 2}\\
\oplus \\
P_2^{\oplus 2}\overset{f_3}{\longrightarrow} P_1^{\oplus 3}
\end{smallmatrix}  \right ]$};
\node[node_style] (v2121) at (8,-2) {$\left [\begin{smallmatrix}
P_2\overset{\mu}{\longrightarrow} P_1\\
\oplus \\
P_2^{\oplus 2}\overset{f_3}{\longrightarrow} P_1^{\oplus 3}
\end{smallmatrix}  \right ]$};
\node[node_style] (v21212) at (10.5,-2) {$\left [\begin{smallmatrix}
P_2\overset{\mu}{\longrightarrow} P_1\\
\oplus \\
P_2\longrightarrow 0
\end{smallmatrix}  \right ]$};
\draw[->]  (v0) edge node{\ } (v1);
\draw[->]  (v1) edge node{\ } (v);
\draw[->]  (v0) edge node{\ } (v2);
\draw[->]  (v2) edge node{\ } (v21);
\draw[->]  (v21) edge node{\ } (v212);
\draw[->]  (v212) edge node{\ } (v2121);
\draw[->]  (v2121) edge node{\ } (v21212);
\draw[->]  (v21212) edge node{\ } (v);
\end{tikzpicture}
\end{center}
where
\begin{center}
$f_1=\left ( \begin{smallmatrix}
\mu\\
\mu\beta\\
\mu\beta^2
\end{smallmatrix} \right )$, $f_2=\left ( \begin{smallmatrix}
\mu\\
\mu\beta
\end{smallmatrix} \right )$, $f_3=\left ( \begin{smallmatrix}
\mu&0\\
-\mu\beta &\mu\\
0&\mu\beta
\end{smallmatrix} \right )$.
\end{center}

Since $\mathsf{Hom}_{\Lambda_4}(P_2,P_1)= e_1\Lambda_4 e_2=K\mu\oplus K\mu\beta\oplus K\mu\beta^2$ and $\mathsf{Hom}_{\Lambda_4}(P_1,P_2)=0$, it is not difficult to compute the left mutations $\mu_{P_1}^-(\Lambda_4)$ and $\mu_{P_2}^-(\Lambda_4)$. According to the bijection introduced in Proposition \ref{tilting-silting}, one can find the corresponding two-term silting complexes. We only show details for the rest of the steps.

(1) Let $T_2=X\oplus Y:=(0\longrightarrow P_1)\oplus(P_2\overset{f_1}{\longrightarrow} P_1^{\oplus 3})$. Then, $\mu_Y^-(T_2)$ does not belong to $\mathsf{2\text{-}silt}\ \Lambda_4$ and therefore, we ignore this mutation. To compute $\mu_X^-(T_2)$, we take a triangle
\begin{center}
$\xymatrix@C=0.7cm{X\ar[r]^{\pi}&Y \ar[r]&\mathsf{cone}(\pi)\ar[r]&X[1]}$ with $\pi=\left ( 0,\left ( \begin{smallmatrix}
0\\
0\\
1
\end{smallmatrix} \right )\right )$.
\end{center}
We may check that $\pi$ is a minimal left $\mathsf{add}(Y)$-approximation. In fact, by the definition,
\begin{itemize}\setlength{\itemsep}{-3pt}
  \item if we compose $\pi$ with the endomorphism
\begin{center}
$\xymatrix@C=1.3cm@R=1cm{Y: & P_2\ar[r]^{f_1}\ar[d]_{k_1e_2+k_2\beta+k_3\beta^2} &P_1^{\oplus 3}\ar[d]^{\left ( \begin{smallmatrix}
k_1&k_2&k_3\\
0&k_1&k_2\\
0&0&k_1
\end{smallmatrix} \right )}\\
Y: & P_2\ar[r]^{f_1} &P_1^{\oplus 3}}$$\xymatrix@R=0.3cm{\ \\ , \text{where }\  k_1,k_2,k_3 \in K,\\\ }$
\end{center}
then all elements of $\mathsf{Hom}_{\mathsf{K^b(proj}\ \Lambda_4)}(X,Y)$ are obtained;
  \item if
$\left ( \begin{smallmatrix}
k_1&k_2&k_3\\
0&k_1&k_2\\
0&0&k_1
\end{smallmatrix} \right )\left ( \begin{smallmatrix}
0\\
0\\
1
\end{smallmatrix} \right )=\left ( \begin{smallmatrix}
0\\
0\\
1
\end{smallmatrix} \right )$, then $k_1=1$ and $k_2=k_3=0$.
\end{itemize}
Hence, $\pi$ is indeed a minimal left $\mathsf{add}(Y)$-approximation. By Lemma \ref{reduced-complex}, we have
\begin{center}
$\mathsf{cone}(\pi)=(\xymatrix@C=1.5cm@R=1.1cm{P_1\oplus P_2\ar[r]^{\ \ \ \left ( \begin{smallmatrix}
0&\mu\\
0&\mu\beta\\
1&\mu\beta^2
\end{smallmatrix} \right )} &P_1^{\oplus 3}})$ $\sim_h(P_2\overset{f_2}{\longrightarrow} P_1^{\oplus 2})$.
\end{center}
Thus, $\mu_X^-(T_2)=(P_2\overset{f_2}{\longrightarrow} P_1^{\oplus 2})\oplus (P_2\overset{f_1}{\longrightarrow} P_1^{\oplus 3})$.

(2) Let $T_{21}=X\oplus Y:=(P_2\overset{f_2}{\longrightarrow} P_1^{\oplus 2})\oplus(P_2\overset{f_1}{\longrightarrow} P_1^{\oplus 3})$. Then, $\mu_X^-(T_{21})\notin \mathsf{2\text{-}silt}\ \Lambda_4$. To compute $\mu_Y^-(T_{21})$, we take a triangle
\begin{center}
$\xymatrix@C=0.7cm{Y\ar[r]^{\pi\ \ }&X^{\oplus3} \ar[r]&\mathsf{cone}(\pi)\ar[r]&Y[1]}$ with $\pi=\left (\left ( \begin{smallmatrix}
e_2\\
\beta\\
\beta^2
\end{smallmatrix} \right ),\left ( \begin{smallmatrix}
1&0&0\\
0&1&0\\
0&1&0\\
0&0&1\\
0&0&1\\
0&0&0
\end{smallmatrix} \right )\right )$.
\end{center}
Then, $\pi$ is a minimal left $\mathsf{add}(X)$-approximation. (In fact, we have $\mathsf{End}_{\mathsf{K^b(proj}\ \Lambda_4)}(X)=K$ since $\mathsf{End}_{\Lambda_4}(P_1)=K$. Thus,
\begin{center}
$\mathsf{End}_{\mathsf{K^b(proj}\ \Lambda_4)}(X^{\oplus3})\simeq \text{Mat}(3,3,K)$.
\end{center}
Secondly, $\lambda \circ \pi =\pi$ for $\lambda\in \text{Mat}(3,3,K)$ implies that $\lambda$ is the identity. Thus, $\pi$ is indeed a minimal left $\mathsf{add}(X)$-approximation.) By applying Lemma \ref{reduced-complex} twice, we have
\begin{center}
$\mathsf{cone}(\pi)=(\xymatrix@C=2.3cm@R=1.1cm{P_2\ar[r]^{\left ( \begin{smallmatrix}
-\mu\\
-\mu\beta\\
-\mu\beta^2\\
e_2\\
\beta\\
\beta^2
\end{smallmatrix} \right )\ \ \ \ \ \  }&P_1^{\oplus3}\oplus P_2^{\oplus3}\ar[r]^{\ \ \ \ \ \ \ \ \left ( \begin{smallmatrix}
1&0&0&\mu&0&0\\
0&1&0&\mu\beta&0&0\\
0&1&0&0&\mu&0\\
0&0&1&0&\mu\beta&0\\
0&0&1&0&0&\mu\\
0&0&0&0&0&\mu\beta
\end{smallmatrix} \right )\ \ \  } &P_1^{\oplus 6}})$ $\sim_h$ $(P_2^{\oplus 2}\overset{f_3}{\longrightarrow} P_1^{\oplus 3})$.
\end{center}
Thus, $\mu_Y^-(T_{21})=(P_2\overset{f_2}{\longrightarrow} P_1^{\oplus 2})\oplus (P_2^{\oplus 2}\overset{f_3}{\longrightarrow} P_1^{\oplus 3})$.

(3) Let $T_{212}=X\oplus Y:=(P_2\overset{f_2}{\longrightarrow} P_1^{\oplus 2})\oplus(P_2^{\oplus 2}\overset{f_3}{\longrightarrow} P_1^{\oplus 3})$. Then, $\mu_Y^-(T_{212})\notin \mathsf{2\text{-}silt}\ \Lambda_4$. To compute $\mu_X^-(T_{212})$, we take a triangle
\begin{center}
$\xymatrix@C=0.7cm{X\ar[r]^{\pi}&Y\ar[r]&\mathsf{cone}(\pi)\ar[r]&X[1]}$ with $\pi=\left ( \left ( \begin{smallmatrix}
0\\
e_2
\end{smallmatrix} \right ),\left ( \begin{smallmatrix}
0&0\\
1&0\\
0&1
\end{smallmatrix} \right )\right )$.
\end{center}
Then, $\pi$ is a minimal left $\mathsf{add}(Y)$-approximation. In fact, if we compose $\pi$ with
\begin{center}
$\xymatrix@C=1.5cm@R=1.2cm{Y: & P_2^{\oplus 2}\ar[r]^{f_3}\ar[d]_{\left ( \begin{smallmatrix}
k_1e_2-k_2\beta&k_2e_2\\
-k_2\beta^2&k_1e_2
\end{smallmatrix} \right )} &P_1^{\oplus 3}\ar[d]^{\left ( \begin{smallmatrix}
k_1&k_2&0\\
0&k_1&-k_2\\
0&0&k_1
\end{smallmatrix} \right )}\\
Y: & P_2^{\oplus 2}\ar[r]^{f_3} &P_1^{\oplus 3}}$$\xymatrix@R=0.3cm{\ \\,\ \text{where }\  k_1,k_2\in K.\\\ }$
\end{center}
then all elements of $\mathsf{Hom}_{\mathsf{K^b(proj}\ \Lambda_4)}(X,Y)$ are obtained; if
$\left ( \begin{smallmatrix}
k_1&k_2&0\\
0&k_1&-k_2\\
0&0&k_1
\end{smallmatrix} \right )\left ( \begin{smallmatrix}
0&0\\
1&0\\
0&1
\end{smallmatrix} \right )=\left ( \begin{smallmatrix}
0&0\\
1&0\\
0&1
\end{smallmatrix} \right )$, then $k_1=1$ and $k_2=0$. By applying Lemma \ref{reduced-complex} twice, we have
\begin{center}
$\mathsf{cone}(\pi)=(\xymatrix@C=2.2cm{P_2\ar[r]^{\left ( \begin{smallmatrix}
-\mu\\
-\mu\beta\\
0\\
e_2
\end{smallmatrix} \right )\ \ \ \ \ \ \ \ \  }&P_1^{\oplus2}\oplus P_2^{\oplus2}\ar[r]^{\ \ \ \ \ \left ( \begin{smallmatrix}
0&0&\mu&0\\
1&0&-\mu\beta&\mu\\
0&1&0&\mu\beta
\end{smallmatrix} \right )} &P_1^{\oplus 3}})$ $\sim_h(P_2\overset{\mu}{\longrightarrow} P_1)$.
\end{center}
Thus, $\mu_X^-(T_{212})=(P_2\overset{\mu}{\longrightarrow} P_1)\oplus (P_2^{\oplus 2}\overset{f_3}{\longrightarrow} P_1^{\oplus 3})$.

(4) Let $T_{2121}=X\oplus Y:=(P_2\overset{\mu}{\longrightarrow} P_1)\oplus(P_2^{\oplus 2}\overset{f_3}{\longrightarrow} P_1^{\oplus 3})$. Then, $\mu_X^-(T_{2121})\notin \mathsf{2\text{-}silt}\ \Lambda_4$. To compute $\mu_Y^-(T_{2121})$, we take a triangle
\begin{center}
$\xymatrix@C=0.7cm{Y\ar[r]^{\pi\ \ }&X^{\oplus3} \ar[r]&\mathsf{cone}(\pi)\ar[r]&Y[1]}$ with $\pi=\left (\left ( \begin{smallmatrix}
e_2&0\\
-\beta&e_2\\
0&\beta
\end{smallmatrix} \right ),\left ( \begin{smallmatrix}
1&0&0\\
0&1&0\\
0&0&1
\end{smallmatrix} \right )\right )$.
\end{center}
Then, $\pi$ is a minimal left $\mathsf{add}(X)$-approximation since $\mathsf{End}_{\mathsf{K^b(proj}\ \Lambda_4)}(X)=K$. Then,
\begin{center}
$\mathsf{cone}(\pi)=(\xymatrix@C=2cm{P_2^{\oplus2}\ar[r]^{\left ( \begin{smallmatrix}
-\mu&0\\
\mu\beta&-\mu\\
0&-\mu\beta\\
e_2&0\\
-\beta&e_2\\
0&\beta
\end{smallmatrix} \right )\ \ \ }&P_1^{\oplus3}\oplus P_2^{\oplus3}\ar[r]^{\ \ \ \ \ \left ( \begin{smallmatrix}
1&0&0&\mu&0&0\\
0&1&0&0&\mu&0\\
0&0&1&0&0&\mu
\end{smallmatrix} \right )} &P_1^{\oplus 3}})$ $\sim_h(P_2\longrightarrow 0)$.
\end{center}
Thus, $\mu_Y^-(T_{2121})=(P_2\overset{\mu}{\longrightarrow} P_1)\oplus (P_2\longrightarrow 0)$.

(5) Let $T_{21212}=X\oplus Y:=(P_2\overset{\mu}{\longrightarrow} P_1)\oplus(P_2\longrightarrow 0)$. Then, it is clear that $\mu_Y^-(T_{21212})$ does not belong to $\mathsf{2\text{-}silt}\ \Lambda_4$ and $\mu_X^-(T_{21212})=(P_1\longrightarrow 0)\oplus (P_2\longrightarrow 0)$.

To sum up the above, we have constructed a finite connected component in $\mathcal{H}(\mathsf{2\text{-}silt}\ \Lambda_4)$. Then, we deduce that $\mathcal{H}(\mathsf{2\text{-}silt}\ \Lambda_4)$ is of type $\mathcal{H}_{1,5}$. By Proposition \ref{tilting-silting}, this is equivalent to saying that $\mathcal{H}(\mathsf{s\tau\text{-}tilt}\ \Lambda_4)$ is of type $\mathcal{H}_{1,5}$.
\end{proof}

\begin{remark}
Let $P_1$ and $P_2$ be the indecomposable projective $\Lambda_3$-modules. We have
\begin{center}
$P_1=\begin{smallmatrix}
e_1\\
\mu\\
\mu\beta
\end{smallmatrix}\simeq\begin{smallmatrix}
1\\
2\\
2
\end{smallmatrix}$ and $P_2=\begin{smallmatrix}
e_2\\
\beta
\end{smallmatrix}\simeq\begin{smallmatrix}
2\\
2
\end{smallmatrix}$.
\end{center}
Then, the Hasse quiver $\mathcal{H}(\mathsf{s\tau\text{-}tilt}\ \Lambda_3)\simeq \mathcal{H}(\mathsf{2\text{-}silt}\ \Lambda_3)$ is of type $\mathcal{H}_{1,3}$ as follows,
\begin{center}
\begin{tikzpicture}[shorten >=1pt, auto, node distance=0cm,
   node_style/.style={font=},
   edge_style/.style={draw=black}]
\node[node_style] (v0) at (0,0) {$\left [\begin{smallmatrix}
0\longrightarrow P_1\\
\oplus \\
0\longrightarrow P_2
\end{smallmatrix}  \right ]$};
\node[node_style] (v1) at (9,0) {$\left [\begin{smallmatrix}
P_1\longrightarrow 0\\
\oplus \\
0\longrightarrow P_2
\end{smallmatrix}  \right ]$};
\node[node_style] (v) at (9,-2) {$\left [\begin{smallmatrix}
P_1\longrightarrow 0\\
\oplus \\
P_2\longrightarrow 0
\end{smallmatrix}  \right ]$};
\node[node_style] (v2) at (0,-2) {$\left [\begin{smallmatrix}
0\longrightarrow P_1\\
\oplus \\
P_2\overset{\binom{\mu}{\mu\beta}}{\longrightarrow} P_1^{\oplus 2}
\end{smallmatrix}  \right ]$};
\node[node_style] (v21) at (3,-2) {$\left [\begin{smallmatrix}
P_2\overset{\mu}{\longrightarrow} P_1\\
\oplus \\
P_2\overset{\binom{\mu}{\mu\beta}}{\longrightarrow} P_1^{\oplus 2}
\end{smallmatrix}  \right ]$};
\node[node_style] (v212) at (6,-2) {$\left [\begin{smallmatrix}
P_2\overset{\mu}{\longrightarrow} P_1\\
\oplus \\
P_2\longrightarrow 0
\end{smallmatrix}  \right ]$};
\draw[->]  (v0) edge node{\ } (v1);
\draw[->]  (v1) edge node{\ } (v);
\draw[->]  (v0) edge node{\ } (v2);
\draw[->]  (v2) edge node{\ } (v21);
\draw[->]  (v21) edge node{\ } (v212);
\draw[->]  (v212) edge node{\ } (v);
\end{tikzpicture}.
\end{center}
\end{remark}

\begin{lemma}
The two-point algebra $\Lambda_5$ is minimal $\tau$-tilting infinite.
\end{lemma}
\begin{proof}
Note that $\Lambda_5$ is a gentle algebra and it is representation-infinite by Hoshino and Miyachi's result \cite[Theorem A]{HM-tame-two-point}. Besides, Plamondon \cite[Theorem 1.1]{P-gentle} showed that a gentle algebra is $\tau$-tilting finite if and only if it is representation-finite. Therefore, $\Lambda_5$ is $\tau$-tilting infinite. For the minimality, we may consider
\begin{center}
$\hat{\Lambda}_5:=\Lambda_5/<\alpha\mu\beta>$
\end{center}
since the socle of $\Lambda_5$ is $K\alpha\mu\beta\oplus K\mu\beta\oplus K\beta$ and any proper quotient $\Lambda_5/I$ of $\Lambda_5$ satisfies $\alpha\mu\beta\in I$. We denote by $P_1$ and $P_2$ the indecomposable projective $\hat{\Lambda}_5$-modules, then
\begin{center}
$P_1=\begin{smallmatrix}
&e_1&\\
\alpha&&\mu\\
\alpha\mu&&\mu\beta
\end{smallmatrix}$ and $P_2=\begin{smallmatrix}
e_2\\
\beta
\end{smallmatrix}$.
\end{center}
By direct calculation, we find that $\mathcal{H}(\mathsf{s\tau\text{-}tilt}\ \hat{\Lambda}_5)\simeq \mathcal{H}(\mathsf{2\text{-}silt}\ \hat{\Lambda}_5))$ is of type $\mathcal{H}_{1,4}$ as follows,
\begin{center}
\begin{tikzpicture}[shorten >=1pt, auto, node distance=0cm,
   node_style/.style={font=},
   edge_style/.style={draw=black}]
\node[node_style] (v0) at (0,0) {$\left [\begin{smallmatrix}
0\longrightarrow P_1\\
\oplus \\
0\longrightarrow P_2
\end{smallmatrix}  \right ]$};
\node[node_style] (v1) at (12,0) {$\left [\begin{smallmatrix}
P_1\longrightarrow 0\\
\oplus \\
0\longrightarrow P_2
\end{smallmatrix}  \right ]$};
\node[node_style] (v) at (12,-2) {$\left [\begin{smallmatrix}
P_1\longrightarrow 0\\
\oplus \\
P_2\longrightarrow 0
\end{smallmatrix}  \right ]$};
\node[node_style] (v2) at (0,-2) {$\left [\begin{smallmatrix}
0\longrightarrow P_1\\
\oplus \\
P_2\overset{\binom{\mu}{\mu\beta}}{\longrightarrow} P_1^{\oplus 2}
\end{smallmatrix}  \right ]$};
\node[node_style] (v21) at (3,-2) {$\left [\begin{smallmatrix}
P_2\overset{\mu}{\longrightarrow} P_1\\
\oplus \\
P_2\overset{\binom{\mu}{\mu\beta}}{\longrightarrow} P_1^{\oplus 2}
\end{smallmatrix}  \right ]$};
\node[node_style] (v212) at (6,-2) {$\left [\begin{smallmatrix}
P_2\overset{\mu}{\longrightarrow} P_1\\
\oplus \\
P_2^{\oplus 2}\overset{(\mu,\alpha\mu)}{\longrightarrow} P_1
\end{smallmatrix}  \right ]$};
\node[node_style] (v2121) at (9,-2) {$\left [\begin{smallmatrix}
P_2\longrightarrow 0\\
\oplus \\
P_2^{\oplus 2}\overset{(\mu,\alpha\mu)}{\longrightarrow} P_1
\end{smallmatrix}  \right ]$};
\draw[->]  (v0) edge node{\ } (v1);
\draw[->]  (v1) edge node{\ } (v);
\draw[->]  (v0) edge node{\ } (v2);
\draw[->]  (v2) edge node{\ } (v21);
\draw[->]  (v21) edge node{\ } (v212);
\draw[->]  (v212) edge node{\ } (v2121);
\draw[->]  (v2121) edge node{\ } (v);
\end{tikzpicture}.
\end{center}
This implies that $\hat{\Lambda}_5$ is $\tau$-tilting finite and $\Lambda_5$ is minimal $\tau$-tilting infinite.
\end{proof}

\begin{lemma}\label{lambda10}
The two-point algebras $\Lambda_7$, $\Lambda_8$, $\Lambda_9$ and $\Lambda_{10}$ are $\tau$-tilting finite.
\end{lemma}
\begin{proof}
Note that $\Lambda_7$, $\Lambda_8$ and $\Lambda_9$ are quotient algebras of $\Lambda_{10}$. By Corollary \ref{factor}, it suffices to show that $\Lambda_{10}$ is $\tau$-tilting finite. The indecomposable projective modules of $\Lambda_{10}$ are
\begin{center}
$P_1=e_1\Lambda_{10}=\substack{ \\\\  \\ \\ \alpha^2\\\alpha^2\mu}\ \substack{\\ \\ \alpha\\ \alpha\mu}\ \substack{e_1\\ \mu \\ \ }\simeq \substack{ \\\\  \\ \\ 1\\2}\ \substack{\\ \\ 1\\2}\ \substack{1\\2 \\ \ }$ and $P_2=e_2\Lambda_{10}=\begin{smallmatrix}
e_2\\
\nu\\
\nu\alpha\\
\nu\alpha^2
\end{smallmatrix}\simeq \begin{smallmatrix}
2\\
1\\
1\\
1
\end{smallmatrix}$.
\end{center}
Since $\mathsf{Hom}_{\Lambda_{10}}(P_1,P_2)=e_2\Lambda_{10} e_1=K\nu\oplus K\nu\alpha\oplus K\nu\alpha^2$ and $\mathsf{Hom}_{\Lambda_{10}}(P_2,P_1)=e_1\Lambda_{10} e_2=K\mu\oplus K\alpha\mu\oplus K\alpha^2\mu$, we know that the computation of the left mutation sequence started at $P_1$ is similar to that of $\Lambda_4$, and the computation of the left mutation sequence started at $P_2$ is similar to that of $\Lambda_4^{\mathsf{op}}$. Then, by Lemma \ref{oppsite algebra} and the calculation in Lemma \ref{lambda4}, we deduce that the Hasse quiver $\mathcal{H}(\mathsf{2\text{-}silt}\ \Lambda_{10})$ is as follows,
\begin{center}
\begin{tikzpicture}[shorten >=1pt, auto, node distance=0cm,
   node_style/.style={font=},
   edge_style/.style={draw=black}]
\node[node_style] (v0) at (0,0) {$\left [\begin{smallmatrix}
0\longrightarrow P_1\\
\oplus \\
0\longrightarrow P_2
\end{smallmatrix}  \right ]$};
\node[node_style] (v1) at (2.5,0) {$\left [\begin{smallmatrix}
P_1\overset{g_1}{\longrightarrow} P_2^{\oplus 3}\\
\oplus \\
0\longrightarrow P_2
\end{smallmatrix}  \right ]$};
\node[node_style] (v12) at (5.2,0) {$\left [\begin{smallmatrix}
P_1\overset{g_1}{\longrightarrow} P_2^{\oplus 3}\\
\oplus \\
P_1\overset{g_2}{\longrightarrow} P_2^{\oplus 2}
\end{smallmatrix}  \right ]$};
\node[node_style] (v121) at (8,0) {$\left [\begin{smallmatrix}
P_1^{\oplus 2}\overset{g_3}{\longrightarrow} P_2^{\oplus 3}\\
\oplus \\
P_1\overset{g_2}{\longrightarrow} P_2^{\oplus 2}
\end{smallmatrix}  \right ]$};
\node[node_style] (v1212) at (11,0) {$\left [\begin{smallmatrix}
P_1^{\oplus 2}\overset{g_3}{\longrightarrow} P_2^{\oplus 3}\\
\oplus \\
P_1\overset{\nu}{\longrightarrow} P_2
\end{smallmatrix}  \right ]$};
\node[node_style] (v12121) at (11,-1.7) {$\left [\begin{smallmatrix}
P_1\longrightarrow 0\\
\oplus \\
P_1\overset{\nu}{\longrightarrow} P_2
\end{smallmatrix}  \right ]$};
\node[node_style] (v) at (11,-3.4) {$\left [\begin{smallmatrix}
P_1\longrightarrow 0\\
\oplus\\
P_2\longrightarrow 0
\end{smallmatrix}  \right ]$};
\node[node_style] (v2) at (0,-1.7) {$\left [\begin{smallmatrix}
0\longrightarrow P_1\\
\oplus \\
P_2\overset{\mu}{\longrightarrow} P_1
\end{smallmatrix}  \right ]$};
\node[node_style] (v21) at (0,-3.4) {$\left [\begin{smallmatrix}
P_2^{\oplus 3}\overset{f_1}{\longrightarrow} P_1^{\oplus 2}\\
\oplus \\
P_2\overset{\mu}{\longrightarrow} P_1
\end{smallmatrix}  \right ]$};
\node[node_style] (v212) at (3,-3.4) {$\left [\begin{smallmatrix}
P_2^{\oplus 3}\overset{f_1}{\longrightarrow} P_1^{\oplus 2}\\
\oplus \\
P_2^{\oplus2}\overset{f_2}{\longrightarrow }P_1
\end{smallmatrix}  \right ]$};
\node[node_style] (v2121) at (5.8,-3.4) {$\left [\begin{smallmatrix}
P_2^{\oplus 3}\overset{f_3}{\longrightarrow} P_1\\
\oplus \\
P_2^{\oplus2}\overset{f_2}{\longrightarrow }P_1
\end{smallmatrix}  \right ]$};
\node[node_style] (v21212) at (8.5,-3.4) {$\left [\begin{smallmatrix}
P_2^{\oplus 3}\overset{f_3}{\longrightarrow} P_1\\
\oplus \\
P_2\longrightarrow 0
\end{smallmatrix}  \right ]$};
\draw[->]  (v0) edge node{\ } (v2);
\draw[->]  (v2) edge node{\ } (v21);
\draw[->]  (v21) edge node{\ } (v212);
\draw[->]  (v212) edge node{\ } (v2121);
\draw[->]  (v2121) edge node{\ } (v21212);
\draw[->]  (v21212) edge node{\ } (v);
\draw[->]  (v0) edge node{\ } (v1);
\draw[->]  (v1) edge node{\ } (v12);
\draw[->]  (v12) edge node{\ } (v121);
\draw[->]  (v121) edge node{\ } (v1212);
\draw[->]  (v1212) edge node{\ } (v12121);
\draw[->]  (v12121) edge node{\ } (v);
\end{tikzpicture},
\end{center}
where $f_1=\left ( \begin{smallmatrix}
\alpha\mu&\mu&0\\
0&-\alpha\mu&\mu
\end{smallmatrix} \right )$, $f_2=\left ( \begin{smallmatrix}
\alpha\mu&\mu
\end{smallmatrix} \right )$, $f_3=\left ( \begin{smallmatrix}
\alpha^2\mu&\alpha\mu&\mu
\end{smallmatrix} \right )$ and
\begin{center}
$g_1=\left ( \begin{smallmatrix}
\nu\\
\nu\alpha\\
\nu\alpha^2
\end{smallmatrix} \right )$, $g_2=\left ( \begin{smallmatrix}
\nu\\
\nu\alpha
\end{smallmatrix} \right )$, $g_3=\left ( \begin{smallmatrix}
\nu&0\\
-\nu\alpha &\nu\\
0&\nu\alpha
\end{smallmatrix} \right )$.
\end{center}

We conclude that $\mathcal{H}(\mathsf{s\tau\text{-}tilt}\ \Lambda_{10})\simeq \mathcal{H}(\mathsf{2\text{-}silt}\ \Lambda_{10})$ is of type $\mathcal{H}_{5,5}$. Thus, $\Lambda_7$, $\Lambda_8$, $\Lambda_9$ and $\Lambda_{10}$ are $\tau$-tilting finite. Next, we determine the type of $\mathcal{H}(\mathsf{s\tau\text{-}tilt}\ \Lambda_i)$ for $i=7,8,9$.

(1) The indecomposable projective $\Lambda_7$-modules are
\begin{center}
$P_1=\begin{smallmatrix}
&e_1&\\
\alpha&&\mu\\
\alpha\mu&&
\end{smallmatrix}$ and $P_2=\begin{smallmatrix}
e_2\\
\nu
\end{smallmatrix}$.
\end{center}
We give the Hasse quiver $\mathcal{H}(\mathsf{2\text{-}silt}\ \Lambda_7)$ by direct calculation as follows,
\begin{center}
\begin{tikzpicture}[shorten >=1pt, auto, node distance=0cm,
   node_style/.style={font=},
   edge_style/.style={draw=black}]
\node[node_style] (v0) at (0,0) {$\left [\begin{smallmatrix}
0\longrightarrow P_1\\
\oplus \\
0\longrightarrow P_2
\end{smallmatrix}  \right ]$};
\node[node_style] (v1) at (3,0) {$\left [\begin{smallmatrix}
P_1\overset{\nu}{\longrightarrow} P_2\\
\oplus \\
0\longrightarrow P_2
\end{smallmatrix}  \right ]$};
\node[node_style] (v12) at (9,0) {$\left [\begin{smallmatrix}
P_1\overset{\nu}{\longrightarrow} P_2\\
\oplus \\
P_1\longrightarrow 0
\end{smallmatrix}  \right ]$};
\node[node_style] (v) at (9,-2) {$\left [\begin{smallmatrix}
P_1\longrightarrow 0\\
\oplus\\
P_2\longrightarrow 0
\end{smallmatrix}  \right ]$};
\node[node_style] (v2) at (0,-2) {$\left [\begin{smallmatrix}
0\longrightarrow P_1\\
\oplus \\
P_2\overset{\mu}{\longrightarrow} P_1
\end{smallmatrix}  \right ]$};
\node[node_style] (v21) at (3,-2) {$\left [\begin{smallmatrix}
P_2^{\oplus 2}\overset{(\mu,\alpha\mu)}{\longrightarrow} P_1\\
\oplus \\
P_2\overset{\mu}{\longrightarrow }P_1
\end{smallmatrix}  \right ]$};
\node[node_style] (v212) at (6,-2) {$\left [\begin{smallmatrix}
P_2^{\oplus 2}\overset{(\mu,\alpha\mu)}{\longrightarrow} P_1\\
\oplus \\
P_2\longrightarrow 0
\end{smallmatrix}  \right ]$};
\draw[->]  (v0) edge node{\ } (v2);
\draw[->]  (v2) edge node{\ } (v21);
\draw[->]  (v21) edge node{\ } (v212);
\draw[->]  (v212) edge node{\ } (v);
\draw[->]  (v0) edge node{\ } (v1);
\draw[->]  (v1) edge node{\ } (v12);
\draw[->]  (v12) edge node{\ } (v);
\end{tikzpicture}.
\end{center}
Hence, $\mathcal{H}(\mathsf{s\tau\text{-}tilt}\ \Lambda_7)\simeq \mathcal{H}(\mathsf{2\text{-}silt}\ \Lambda_{7})$ is of type $\mathcal{H}_{2,3}$.

(2) The indecomposable projective $\Lambda_8$-modules are
\begin{center}
$P_1=\begin{smallmatrix}
&e_1&\\
\alpha&&\mu\\
\alpha\mu&&
\end{smallmatrix}$ and $P_2=\begin{smallmatrix}
e_2\\
\nu\\
\nu\alpha
\end{smallmatrix}$,
\end{center}
and the Hasse quiver $\mathcal{H}(\mathsf{2\text{-}silt}\ \Lambda_8)$ is given as follows,
\begin{center}
\begin{tikzpicture}[shorten >=1pt, auto, node distance=0cm,
   node_style/.style={font=},
   edge_style/.style={draw=black}]
\node[node_style] (v0) at (0,0) {$\left [\begin{smallmatrix}
0\longrightarrow P_1\\
\oplus \\
0\longrightarrow P_2
\end{smallmatrix}  \right ]$};
\node[node_style] (v1) at (3,0) {$\left [\begin{smallmatrix}
P_1\overset{\binom{\nu}{\nu\alpha}}{\longrightarrow} P_2^{\oplus 2}\\
\oplus \\
0\longrightarrow P_2
\end{smallmatrix}  \right ]$};
\node[node_style] (v12) at (6,0) {$\left [\begin{smallmatrix}
P_1\overset{\binom{\nu}{\nu\alpha}}{\longrightarrow} P_2^{\oplus 2}\\
\oplus \\
P_1\overset{\nu}{\longrightarrow} P_2\\
\end{smallmatrix}  \right ]$};
\node[node_style] (v121) at (9,0) {$\left [\begin{smallmatrix}
P_1\longrightarrow 0\\
\oplus \\
P_1\overset{\nu}{\longrightarrow} P_2
\end{smallmatrix}  \right ]$};
\node[node_style] (v) at (9,-2) {$\left [\begin{smallmatrix}
P_1\longrightarrow 0\\
\oplus\\
P_2\longrightarrow 0
\end{smallmatrix}  \right ]$};
\node[node_style] (v2) at (0,-2) {$\left [\begin{smallmatrix}
0\longrightarrow P_1\\
\oplus \\
P_2\overset{\mu}{\longrightarrow} P_1
\end{smallmatrix}  \right ]$};
\node[node_style] (v21) at (3,-2) {$\left [\begin{smallmatrix}
P_2^{\oplus 2}\overset{(\mu,\alpha\mu)}{\longrightarrow} P_1\\
\oplus \\
P_2\overset{\mu}{\longrightarrow }P_1
\end{smallmatrix}  \right ]$};
\node[node_style] (v212) at (6,-2) {$\left [\begin{smallmatrix}
P_2^{\oplus 2}\overset{(\mu,\alpha\mu)}{\longrightarrow} P_1\\
\oplus \\
P_2\longrightarrow 0
\end{smallmatrix}  \right ]$};
\draw[->]  (v0) edge node{\ } (v2);
\draw[->]  (v2) edge node{\ } (v21);
\draw[->]  (v21) edge node{\ } (v212);
\draw[->]  (v212) edge node{\ } (v);
\draw[->]  (v0) edge node{\ } (v1);
\draw[->]  (v1) edge node{\ } (v12);
\draw[->]  (v12) edge node{\ } (v121);
\draw[->]  (v121) edge node{\ } (v);
\end{tikzpicture}.
\end{center}
Hence, $\mathcal{H}(\mathsf{s\tau\text{-}tilt}\ \Lambda_8)\simeq \mathcal{H}(\mathsf{2\text{-}silt}\ \Lambda_{8})$ is of type $\mathcal{H}_{3,3}$.

(3) Let $Q_1$ and $Q_2$ be the indecomposable projective $\Lambda_9$-modules. Then,
\begin{center}
$Q_1=e_1\Lambda_9=\substack{ \\\\  \\ \\ \alpha^2\\\alpha^2\mu}\ \substack{\\ \\ \alpha\\ \alpha\mu}\ \substack{e_1\\ \mu \\ \ }\simeq \substack{ \\\\  \\ \\ 1\\2}\ \substack{\\ \\ 1\\2}\ \substack{1\\2 \\ \ }$ and $Q_2=e_2\Lambda_9=\begin{smallmatrix}
e_2\\
\nu
\end{smallmatrix}\simeq \begin{smallmatrix}
2\\
1
\end{smallmatrix}$.
\end{center}
Since $\mathsf{Hom}_{\Lambda_9}(Q_1,Q_2)= e_2\Lambda_9 e_1=K\nu$ and $\mathsf{Hom}_{\Lambda_9}(Q_2,Q_1)=e_1\Lambda_9 e_2=K\mu\oplus K\alpha\mu\oplus K\alpha^2\mu$, the computation of the left mutation sequence started at $Q_2$ is similar to that of $\Lambda_4^{\mathsf{op}}$. Then, by Lemma \ref{oppsite algebra} and the calculation in Lemma \ref{lambda4}, we deduce that
$\mathcal{H}(\mathsf{2\text{-}silt}\ \Lambda_9)$ is presented as follows,
\begin{center}
\begin{tikzpicture}[shorten >=1pt, auto, node distance=0cm,
   node_style/.style={font=},
   edge_style/.style={draw=black}]
\node[node_style] (v0) at (0,0) {$\left [\begin{smallmatrix}
0\longrightarrow Q_1\\
\oplus \\
0\longrightarrow Q_2
\end{smallmatrix}  \right ]$};
\node[node_style] (v1) at (5,0) {$\left [\begin{smallmatrix}
Q_1\overset{\nu}{\longrightarrow} Q_2\\
\oplus \\
0\longrightarrow Q_2
\end{smallmatrix}  \right ]$};
\node[node_style] (v12) at (11,0) {$\left [\begin{smallmatrix}
Q_1\longrightarrow 0\\
\oplus \\
Q_1\overset{\nu}{\longrightarrow} Q_2
\end{smallmatrix}  \right ]$};
\node[node_style] (v) at (11,-4) {$\left [\begin{smallmatrix}
Q_1\longrightarrow 0\\
\oplus\\
Q_2\longrightarrow 0
\end{smallmatrix}  \right ]$};
\node[node_style] (v2) at (0,-2) {$\left [\begin{smallmatrix}
0\longrightarrow Q_1\\
\oplus \\
Q_2\overset{\mu}{\longrightarrow} Q_1
\end{smallmatrix}  \right ]$};
\node[node_style] (v21) at (0,-4) {$\left [\begin{smallmatrix}
Q_2^{\oplus 3}\overset{f_1}{\longrightarrow} Q_1^{\oplus 2}\\
\oplus \\
Q_2\overset{\mu}{\longrightarrow} Q_1
\end{smallmatrix}  \right ]$};
\node[node_style] (v212) at (3,-4) {$\left [\begin{smallmatrix}
Q_2^{\oplus 3}\overset{f_1}{\longrightarrow} Q_1^{\oplus 2}\\
\oplus \\
Q_2^{\oplus2}\overset{f_2}{\longrightarrow }Q_1
\end{smallmatrix}  \right ]$};
\node[node_style] (v2121) at (5.8,-4) {$\left [\begin{smallmatrix}
Q_2^{\oplus 3}\overset{f_3}{\longrightarrow} Q_1\\
\oplus \\
Q_2^{\oplus2}\overset{f_2}{\longrightarrow }Q_1
\end{smallmatrix}  \right ]$};
\node[node_style] (v21212) at (8.5,-4) {$\left [\begin{smallmatrix}
Q_2^{\oplus 3}\overset{f_3}{\longrightarrow} Q_1\\
\oplus \\
Q_2\longrightarrow 0
\end{smallmatrix}  \right ]$};
\draw[->]  (v0) edge node{\ } (v2);
\draw[->]  (v2) edge node{\ } (v21);
\draw[->]  (v21) edge node{\ } (v212);
\draw[->]  (v212) edge node{\ } (v2121);
\draw[->]  (v2121) edge node{\ } (v21212);
\draw[->]  (v21212) edge node{\ } (v);
\draw[->]  (v0) edge node{\ } (v1);
\draw[->]  (v1) edge node{\ } (v12);
\draw[->]  (v12) edge node{\ } (v);
\end{tikzpicture},
\end{center}
where $f_1=\left ( \begin{smallmatrix}
\alpha\mu&\mu&0\\
0&-\alpha\mu&\mu
\end{smallmatrix} \right )$, $f_2=\left ( \begin{smallmatrix}
\alpha\mu&\mu
\end{smallmatrix} \right )$ and $f_3=\left ( \begin{smallmatrix}
\alpha^2\mu&\alpha\mu&\mu
\end{smallmatrix} \right )$. By Proposition \ref{tilting-silting}, we conclude that $\mathcal{H}(\mathsf{s\tau\text{-}tilt}\ \Lambda_9)\simeq \mathcal{H}(\mathsf{2\text{-}silt}\ \Lambda_9)$ is of type $\mathcal{H}_{2,5}$.
\end{proof}

\begin{lemma}
The two-point algebras $\Lambda_{11}$ and $\Lambda_{12}$ are $\tau$-tilting finite.
\end{lemma}
\begin{proof}
(1) The indecomposable projective $\Lambda_{11}$-modules are
\begin{center}
$P_1=\begin{smallmatrix}
&e_1&\\
\alpha&&\mu\\
&&\mu\beta
\end{smallmatrix}$ and $P_2=\begin{smallmatrix}
&e_2&\\
\beta&&\nu\\
&&\nu\alpha
\end{smallmatrix}$.
\end{center}
We calculate the Hasse quiver $\mathcal{H}(\mathsf{2\text{-}silt}\ \Lambda_{11})$ directly as follows,
\begin{center}
\begin{tikzpicture}[shorten >=1pt, auto, node distance=0cm,
   node_style/.style={font=},
   edge_style/.style={draw=black}]
\node[node_style] (v0) at (0,0) {$\left [\begin{smallmatrix}
0\longrightarrow P_1\\
\oplus \\
0\longrightarrow P_2
\end{smallmatrix}  \right ]$};
\node[node_style] (v1) at (3,0) {$\left [\begin{smallmatrix}
P_1\overset{\binom{\nu}{\nu\alpha}}{\longrightarrow} P_2^{\oplus 2}\\
\oplus \\
0\longrightarrow P_2
\end{smallmatrix}  \right ]$};
\node[node_style] (v12) at (6,0) {$\left [\begin{smallmatrix}
P_1\overset{\binom{\nu}{\nu\alpha}}{\longrightarrow} P_2^{\oplus 2}\\
\oplus \\
P_1\overset{\nu}{\longrightarrow} P_2\\
\end{smallmatrix}  \right ]$};
\node[node_style] (v121) at (9,0) {$\left [\begin{smallmatrix}
P_1\longrightarrow 0\\
\oplus \\
P_1\overset{\nu}{\longrightarrow} P_2
\end{smallmatrix}  \right ]$};
\node[node_style] (v) at (9,-2) {$\left [\begin{smallmatrix}
P_1\longrightarrow 0\\
\oplus \\
P_2\longrightarrow 0
\end{smallmatrix}  \right ]$};
\node[node_style] (v2) at (0,-2) {$\left [\begin{smallmatrix}
0\longrightarrow P_1\\
\oplus \\
P_2\overset{\binom{\mu}{\mu\beta}}{\longrightarrow} P_1^{\oplus 2}
\end{smallmatrix}  \right ]$};
\node[node_style] (v21) at (3,-2) {$\left [\begin{smallmatrix}
P_2\overset{\mu}{\longrightarrow} P_1\\
\oplus \\
P_2\overset{\binom{\mu}{\mu\beta}}{\longrightarrow} P_1^{\oplus 2}
\end{smallmatrix}  \right ]$};
\node[node_style] (v212) at (6,-2) {$\left [\begin{smallmatrix}
P_2\overset{\mu}{\longrightarrow} P_1\\
\oplus \\
P_2\longrightarrow 0
\end{smallmatrix}  \right ]$};
\draw[->]  (v0) edge node{\ } (v1);
\draw[->]  (v1) edge node{\ } (v12);
\draw[->]  (v12) edge node{\ } (v121);
\draw[->]  (v121) edge node{\ } (v);
\draw[->]  (v0) edge node{\ } (v2);
\draw[->]  (v2) edge node{\ } (v21);
\draw[->]  (v21) edge node{\ } (v212);
\draw[->]  (v212) edge node{\ } (v);
\end{tikzpicture}.
\end{center}
Thus, $\mathcal{H}(\mathsf{s\tau\text{-}tilt}\ \Lambda_{11})\simeq \mathcal{H}(\mathsf{2\text{-}silt}\ \Lambda_{11})$ is of type $\mathcal{H}_{3,3}$.

(2) The indecomposable projective $\Lambda_{12}$-modules are
\begin{center}
$P_1=\begin{smallmatrix}
&e_1&\\
\alpha&&\mu\\
\alpha\mu&&\mu\beta
\end{smallmatrix}$ and $P_2=\begin{smallmatrix}
&e_2&\\
\beta&&\nu
\end{smallmatrix}$,
\end{center}
and the Hasse quiver $\mathcal{H}(\mathsf{2\text{-}silt}\ \Lambda_{12})$ is shown as follows,
\begin{center}
\begin{tikzpicture}[shorten >=1pt, auto, node distance=0cm,
   node_style/.style={font=},
   edge_style/.style={draw=black}]
\node[node_style] (v0) at (0,0) {$\left [\begin{smallmatrix}
0\longrightarrow P_1\\
\oplus \\
0\longrightarrow P_2
\end{smallmatrix}  \right ]$};
\node[node_style] (v1) at (6,0) {$\left [\begin{smallmatrix}
P_1\overset{\nu}{\longrightarrow} P_2\\
\oplus \\
0\longrightarrow P_2
\end{smallmatrix}  \right ]$};
\node[node_style] (v12) at (12,0) {$\left [\begin{smallmatrix}
P_1\overset{\nu}{\longrightarrow} P_2\\
\oplus \\
P_1\longrightarrow 0
\end{smallmatrix}  \right ]$};
\node[node_style] (v) at (12,-2) {$\left [\begin{smallmatrix}
P_1\longrightarrow 0\\
\oplus \\
P_2\longrightarrow 0
\end{smallmatrix}  \right ]$};
\node[node_style] (v2) at (0,-2) {$\left [\begin{smallmatrix}
0\longrightarrow P_1\\
\oplus \\
P_2\overset{\binom{\mu}{\mu\beta}}{\longrightarrow} P_1^{\oplus 2}
\end{smallmatrix}  \right ]$};
\node[node_style] (v21) at (3,-2) {$\left [\begin{smallmatrix}
P_2\overset{\mu}{\longrightarrow} P_1\\
\oplus \\
P_2\overset{\binom{\mu}{\mu\beta}}{\longrightarrow} P_1^{\oplus 2}
\end{smallmatrix}  \right ]$};
\node[node_style] (v212) at (6,-2) {$\left [\begin{smallmatrix}
P_2\overset{\mu}{\longrightarrow} P_1\\
\oplus \\
P_2^{\oplus 2}\overset{(\mu,\alpha\mu)}{\longrightarrow} P_1
\end{smallmatrix}  \right ]$};
\node[node_style] (v2121) at (9,-2) {$\left [\begin{smallmatrix}
P_2\longrightarrow 0\\
\oplus \\
P_2^{\oplus 2}\overset{(\mu,\alpha\mu)}{\longrightarrow} P_1
\end{smallmatrix}  \right ]$};
\draw[->]  (v0) edge node{\ } (v1);
\draw[->]  (v1) edge node{\ } (v12);
\draw[->]  (v12) edge node{\ } (v);
\draw[->]  (v0) edge node{\ } (v2);
\draw[->]  (v2) edge node{\ } (v21);
\draw[->]  (v21) edge node{\ } (v212);
\draw[->]  (v212) edge node{\ } (v2121);
\draw[->]  (v2121) edge node{\ } (v);
\end{tikzpicture}.
\end{center}
Thus, $\mathcal{H}(\mathsf{s\tau\text{-}tilt}\ \Lambda_{12})\simeq \mathcal{H}(\mathsf{2\text{-}silt}\ \Lambda_{12})$ is of type $\mathcal{H}_{2,4}$.
\end{proof}

Lastly, we summarize the number $\#\mathsf{s\tau\text{-}tilt}\ \Lambda_i$ and the type of $\mathcal{H}(\mathsf{s\tau\text{-}tilt}\ \Lambda_i)$ for $\Lambda_i$ with $i\neq 2, 5$ as follows,
\begin{center}
\renewcommand\arraystretch{1.2}
\begin{tabular}{|c|c|c|c|c|c|c|c|c|c|c|c|}
\hline
$\Lambda_i$&$\Lambda_1$&$\Lambda_3$&$\Lambda_4$&$\Lambda_6$&$\Lambda_7$&$\Lambda_8$&$\Lambda_9$&$\Lambda_{10}$&$\Lambda_{11}$& $\Lambda_{12}$\\ \hline
$\#\mathsf{s\tau\text{-}tilt}\ \Lambda_i$&5&6& 8&6&7&8 &9 & 12&8& 8\\ \hline
Type&$\mathcal{H}_{1,2}$ &$\mathcal{H}_{1,3}$&$\mathcal{H}_{1,5}$&$\mathcal{H}_{2,2}$&$\mathcal{H}_{2,3}$&$\mathcal{H}_{3,3}$&$\mathcal{H}_{2,5}$&$\mathcal{H}_{5,5}$&$\mathcal{H}_{3,3}$&$\mathcal{H}_{2,4}$ \\ \hline
\end{tabular}.
\end{center}

\subsection{The proof of Theorem \ref{table-W}}
First, one can easily find that $W_1$, $W_2$, $W_3$ and $W_5$ have $\Lambda_2$ as a quotient algebra and therefore, they are $\tau$-tilting infinite. It is also not difficult to find that $W_4$ is $\tau$-tilting finite and $\mathcal{H}(\mathsf{s\tau\text{-}tilt}\ W_4)$ is of type $\mathcal{H}_{1,2}$.

We may distinguish the case $W_{14}$. Note that $\nu\alpha\mu\in C(W_{14})$ and the indecomposable projective modules of $\widetilde{W}_{14}:=W_{14}/<\nu\alpha\mu>$ are
\begin{center}
$P_1=\begin{smallmatrix}
&&e_1&\\
&\alpha&&\mu\\
\alpha^2&&\alpha\mu&
\end{smallmatrix}$ and $P_2=\begin{smallmatrix}
e_2\\
\nu\\
\nu\alpha\\
\nu\alpha^2
\end{smallmatrix}$.
\end{center}
Then, we find that $\widetilde{W}_{14}$ is a quotient algebra of $\Lambda_{10}$ by $\alpha^2\mu$. Thus, by similar calculation with $\Lambda_{10}$ in the proof of Lemma \ref{lambda10}, one can check that $\mathcal{H}(\mathsf{s\tau\text{-}tilt}\ \widetilde{W}_{14})$ is of type $\mathcal{H}_{3,5}$. By Lemma \ref{centers}, we deduce that $\mathcal{H}(\mathsf{s\tau\text{-}tilt}\ W_{14})\simeq \mathcal{H}(\mathsf{s\tau\text{-}tilt}\ \widetilde{W}_{14})$ is of type $\mathcal{H}_{3,5}$.

Next, we show that $W_6\sim W_{34}$ (except for $W_{14}$) are $\tau$-tilting finite by determining the type of $\mathcal{H}(\mathsf{s\tau\text{-}tilt}\ W_i)$ for $i=6, 7, \ldots, 34$ ($i\neq 14$). In order to do this, we can apply Lemma \ref{centers} to construct a two-sided ideal $I$ generated by elements in $C(W_i)\cap \mathsf{rad}(W_i)$ such that $\mathsf{s\tau\text{-}tilt}\ A\simeq \mathsf{s\tau\text{-}tilt}\ \left ( W_i/I \right )$. Then, we can find the type of $\mathcal{H}(\mathsf{s\tau\text{-}tilt}\ W_i)$ following Table $\Lambda$. Here, we compute the center of an algebra by GAP as follows, see \cite{gap}.
\begin{center}
\renewcommand\arraystretch{1.1}
\begin{tabular}{|c|c|c|}
\hline
$W_i$&$I$&$A$  \\ \hline
$W_6$&$\alpha^2$&$\Lambda_3^{\mathsf{op}}$ \\ \hline
$W_7$&$\alpha^3$&$\Lambda_4^{\mathsf{op}}$ \\ \hline
$W_8$&$\alpha$&$\Lambda_4$ \\ \hline
$W_9$&\multirow{2}*{$\alpha,\beta^2$}&\multirow{2}*{$\Lambda_3$}  \\
\cline{1-1}
$W_{10}$&& \\ \hline
$W_{11}$&$\beta^2$&$\hat{\Lambda}_5$ \\ \hline
$W_{12}$&\multirow{2}*{$\alpha,\beta$}&\multirow{2}*{$\Lambda_1$}  \\
\cline{1-1}
$W_{13}$&& \\ \hline
$W_{15}$&$\nu\mu$&$\Lambda_9^{\mathsf{op}}$ \\ \hline
$W_{16}$&$\alpha^2,\nu\mu$&$\Lambda_8$ \\ \hline
$W_{17}$&$\alpha^3$&$\Lambda_9^{\mathsf{op}}$ \\ \hline
$W_{18}$&\multirow{2}*{$\alpha^2$}&$\Lambda_8$ \\
\cline{1-1}\cline{3-3}
$W_{19}$&&\multirow{2}*{$\Lambda_7$} \\
\cline{1-2}
$W_{20}$&$\mu\nu+\nu\mu$& \\ \hline
\end{tabular}\ \ \ \
\renewcommand\arraystretch{1.1}
\begin{tabular}{|c|c|c|}
\hline
$W_i$&$I$&$A$  \\ \hline
$W_{21}$&$\alpha\mu\nu$&\multirow{2}*{$\Lambda_7$}\\
\cline{1-2}
$W_{22}$&$\alpha^2,\mu\nu$& \\ \hline
$W_{23}$&$\alpha^2+\nu\mu,\nu\alpha\mu$&$\Lambda_8$ \\ \hline
$W_{24}$&$\mu\nu$&$\widetilde{W}_{14}^{\mathsf{op}}$ \\ \hline
$W_{25}$&$\alpha^2,\beta$&$\Lambda_7$ \\ \hline
$W_{26}$&$\alpha,\mu\nu$&$\Lambda_7^{\mathsf{op}}$ \\ \hline
$W_{27}$&\multirow{4}*{$\mu\nu$}&$\Lambda_{12}^{\mathsf{op}}$ \\
\cline{1-1}\cline{3-3}
$W_{28}$&&$\Lambda_{11}$\\
\cline{1-1}\cline{3-3}
$W_{29}$&&$\Lambda_{11}^{\mathsf{op}}$ \\
\cline{1-1}\cline{3-3}
$W_{30}$&&$\Lambda_{12}$ \\ \hline
$W_{31}$&$\alpha+\beta,\nu\mu$&\multirow{4}*{$\Lambda_6$} \\
\cline{1-2}
$W_{32}$&$\alpha+\beta,\nu\mu,\mu\nu$& \\
\cline{1-2}
$W_{33}$&$\alpha+\beta$& \\
\cline{1-2}
$W_{34}$&$\alpha+\beta,\mu\nu$& \\ \hline
\end{tabular}
\end{center}
In particular, we point out that although $\Lambda_7\not\simeq \widetilde{W}_{21}:=W_{21}/I$, but $\mathsf{s\tau\text{-}tilt}\ \Lambda_7\simeq \mathsf{s\tau\text{-}tilt}\ \widetilde{W}_{21}$. To see the latter one, one may check that $\nu\mu+\mu\nu \in C(\widetilde{W}_{21})$ and therefore,
\begin{center}
$\mathsf{s\tau\text{-}tilt}\ \widetilde{W}_{21}\simeq \mathsf{s\tau\text{-}tilt}\ \left (\widetilde{W}_{21}/<\mu\nu,\nu\mu>\right )\simeq \mathsf{s\tau\text{-}tilt}\ \Lambda_7$.
\end{center}

\subsection{The proof of Theorem \ref{table-T}}
Similar to the proof of Theorem \ref{table-W}, one can check that $T_1$ has $\Lambda_2$ as a quotient algebra, $T_3$ and $T_{17}$ have $\Lambda_5$ as a quotient algebra. Hence, $T_1$, $T_3$ and $T_{17}$ are $\tau$-tilting infinite. We may also distinguish the following cases.

Case ($T_9$). Since $\nu\mu, \alpha\mu\nu+\mu\nu\alpha+\nu\alpha\mu \in C(T_9)$ and $\alpha\mu\nu\in C(\widetilde{T}_9)$ with
\begin{center}
$\widetilde{T}_9:=T_9/<\nu\mu,\nu\alpha\mu,\alpha\mu\nu+\mu\nu\alpha>$.
\end{center}
Then, we have $\mu\nu\in C(\widetilde{T}_9/<\alpha\mu\nu>)$ and therefore,
\begin{center}
$\mathsf{s\tau\text{-}tilt}\ T_9\simeq \mathsf{s\tau\text{-}tilt}\ \left (T_9/<\mu\nu,\nu\mu,\nu\alpha\mu>\right )\simeq \mathsf{s\tau\text{-}tilt}\ \Lambda_8$.
\end{center}

Case ($T_{20}$). For any $k\in K/\{0\}$, we have $\mu\nu+\nu\mu \in C(T_{20})$ such that
\begin{center}
$\mathsf{s\tau\text{-}tilt}\ T_{20}\simeq \mathsf{s\tau\text{-}tilt}\ \widetilde{T}_{20}$ with $\widetilde{T}_{20}:=T_{20}/<\mu\nu,\nu\mu>$.
\end{center}
Then, the indecomposable projective $\widetilde{T}_{20}$-modules are
\begin{center}
$P_1=\begin{smallmatrix}
&e_1&\\
\alpha&&\mu\\
&\mu\beta&
\end{smallmatrix}$ and $P_2=\begin{smallmatrix}
&e_2&\\
\beta&&\nu\\
&&\nu\alpha
\end{smallmatrix}$,
\end{center}
and the Hasse quiver $\mathcal{H}(\mathsf{2\text{-}silt}\ \widetilde{T}_{20})$ is given as follows,
\begin{center}
\begin{tikzpicture}[shorten >=1pt, auto, node distance=0cm,
   node_style/.style={font=},
   edge_style/.style={draw=black}]
\node[node_style] (v0) at (0,0) {$\left [\begin{smallmatrix}
0\longrightarrow P_1\\
\oplus \\
0\longrightarrow P_2
\end{smallmatrix}  \right ]$};
\node[node_style] (v1) at (3,0) {$\left [\begin{smallmatrix}
P_1\overset{\binom{\nu}{\nu\alpha}}{\longrightarrow} P_2^{\oplus 2}\\
\oplus \\
0\longrightarrow P_2
\end{smallmatrix}  \right ]$};
\node[node_style] (v12) at (6,0) {$\left [\begin{smallmatrix}
P_1\overset{\binom{\nu}{\nu\alpha}}{\longrightarrow} P_2^{\oplus 2}\\
\oplus \\
P_1\overset{\nu}{\longrightarrow} P_2\\
\end{smallmatrix}  \right ]$};
\node[node_style] (v121) at (9,0) {$\left [\begin{smallmatrix}
P_1\longrightarrow 0\\
\oplus \\
P_1\overset{\nu}{\longrightarrow} P_2
\end{smallmatrix}  \right ]$};
\node[node_style] (v) at (9,-2) {$\left [\begin{smallmatrix}
P_1\longrightarrow 0\\
\oplus \\
P_2\longrightarrow 0
\end{smallmatrix}  \right ]$};
\node[node_style] (v2) at (0,-2) {$\left [\begin{smallmatrix}
0\longrightarrow P_1\\
\oplus \\
P_2\overset{\mu}{\longrightarrow} P_1
\end{smallmatrix}  \right ]$};
\node[node_style] (v21) at (6,-2) {$\left [\begin{smallmatrix}
P_2\overset{\mu}{\longrightarrow} P_1\\
\oplus \\
P_2\longrightarrow 0
\end{smallmatrix}  \right ]$};
\draw[->]  (v0) edge node{\ } (v1);
\draw[->]  (v1) edge node{\ } (v12);
\draw[->]  (v12) edge node{\ } (v121);
\draw[->]  (v121) edge node{\ } (v);
\draw[->]  (v0) edge node{\ } (v2);
\draw[->]  (v2) edge node{\ } (v21);
\draw[->]  (v21) edge node{\ } (v);
\end{tikzpicture}.
\end{center}
Thus, $\mathcal{H}(\mathsf{s\tau\text{-}tilt}\ T_{20})\simeq \mathcal{H}(\mathsf{s\tau\text{-}tilt}\ \widetilde{T}_{20})\simeq \mathcal{H}(\mathsf{2\text{-}silt}\ \widetilde{T}_{20})$ is of type $\mathcal{H}_{2,3}$.

Case ($T_{21}$). For any $k_1,k_2\in K/\{0\}$, we have $\mu\nu+\nu\mu \in C(T_{21})$. Similarly, we have $\mathsf{s\tau\text{-}tilt}\ T_{21}\simeq \mathsf{s\tau\text{-}tilt}\ \left (T_{21}/<\mu\nu,\nu\mu>\right )$ and the corresponding Hasse quiver is of type $\mathcal{H}_{2,2}$.

For the remaining cases, we may apply Lemma \ref{centers} to construct a two-sided ideal $I$ generated by elements in $C(T_i)\cap \mathsf{rad}(T_i)$ such that $\mathsf{s\tau\text{-}tilt}\ B\simeq \mathsf{s\tau\text{-}tilt}\ \left ( T_i/I \right )$, as follows,
\begin{center}
\renewcommand\arraystretch{1.2}
\begin{tabular}{|c|c|c|}
\hline
$T_i$&$I$&$B$  \\ \hline
$T_2$&$\alpha^2$&\multirow{2}*{$\Lambda_3^{\mathsf{op}}$}  \\
\cline{1-2}
$T_4$&$\beta$&   \\ \hline
$T_5$&$\alpha, \beta$&$\Lambda_1$  \\ \hline
$T_6$&$\alpha+\beta^2$&$\Lambda_3$ \\ \hline
$T_7$&\multirow{2}*{$\alpha+\beta$}&\multirow{2}*{$\Lambda_1$}  \\
\cline{1-1}
$T_8$&& \\ \hline
$T_{10}$&$\nu\alpha^2\mu$&$\Lambda_{10}$ \\ \hline
$T_{11}$&$\alpha^2,\nu\alpha\mu$&$\Lambda_8$ \\ \hline

\end{tabular}
\renewcommand\arraystretch{1.2}
\begin{tabular}{|c|c|c|}
\hline
$T_i$&$I$&$B$  \\ \hline
$T_{12}$&$\alpha^2,\nu\mu$&$\Lambda_7$ \\ \hline
$T_{13}$&$\alpha,\mu\nu+\nu\mu$&$\Lambda_6$ \\ \hline
$T_{14}$&$\alpha^2+\nu\mu$&$\Lambda_8$ \\ \hline
$T_{15}$&$\alpha^2, \nu\mu$&$\Lambda_7$ \\ \hline
$T_{16}$&$\mu\nu$&$\Lambda_{9}$ \\ \hline
\multirow{2}*{$T_{18}$}&$\beta,\nu\alpha\mu+$&\multirow{2}*{$\widetilde{T}_9$} \\
&$\alpha\mu\nu+\mu\nu\alpha$& \\\hline
$T_{19}$&$\alpha,\beta,\mu\nu+\nu\mu$&$\Lambda_6$ \\ \hline
\end{tabular}.
\end{center}

\subsection{Other applications}
At the end of this paper, we give two easy observations. First, we have
\begin{proposition}
Let $\Lambda$ be a connected two-point algebra without loops. Then, $\Lambda$ is $\tau$-tilting finite if and only if it is representation-finite.
\end{proposition}
\begin{proof}
By our assumption, the quiver $Q$ of $\Lambda$ does not contain loops. If $Q$ contains multiple arrows, then $\Lambda$ has the Kronecker algebra $\Lambda_2$ as a quotient algebra and hence, $\Lambda$ is $\tau$-tilting infinite. Then, we deduce that if $\Lambda$ is $\tau$-tilting finite, then $Q$ is either $\xymatrix@C=0.7cm{\bullet\ar@<0.5ex>[r]^{\ }&\bullet\ar@<0.5ex>[l]_{\ }}$ or $\xymatrix@C=0.7cm{\bullet\ar[r]&\bullet}$. On the other hand, any finite-dimensional algebra with quiver $\xymatrix@C=0.7cm{\bullet\ar@<0.5ex>[r]^{\ }&\bullet\ar@<0.5ex>[l]_{\ }}$ or $\xymatrix@C=0.7cm{\bullet\ar[r]&\bullet}$ is representation-finite from Bongartz and Gabriel \cite{BG-covering}.
\end{proof}

Second, we determine the $\tau$-tilting finiteness of two-point symmetric special biserial algebras. We refer to \cite{S-brauer graph} for the basic concepts and properties of symmetric special biserial algebras, or equivalently, Brauer graph algebras. In \cite{AIP-joint}, the authors classified two-point symmetric special biserial algebras up to Morita equivalence, so that we can determine their $\tau$-tilting finiteness.
\begin{proposition}{\rm{(\cite[Theorem 7.1]{AIP-joint})}}
Let $\Lambda$ be a two-point symmetric special biserial algebra. Then, $\Lambda$ is Morita equivalent to one of $B_i=KQ/I_i$ below.

\ \ \ \ \ $Q: \xymatrix{\bullet \ar@<0.5ex>[r]^{\mu}&\bullet\ar@<0.5ex>[l]^{\nu}}$
\vspace{0.3cm}

$I_1:$ $(\mu\nu)^n\mu=(\nu\mu)^n\nu=0, n\geqslant 1$.

\rule[0.5\baselineskip]{6cm}{1pt}

\ \ \ \ \ $Q: \xymatrix@C=1cm{\bullet \ar@<0.3ex>[r]\ar@<1ex>[r]^{\mu_1,\mu_2}&\bullet\ar@<0.3ex>[l]\ar@<1ex>[l]^{\nu_1,\nu_2}}$
\vspace{0.3cm}

$I_2:$ $\mu_1\nu_2=\nu_2\mu_1=\mu_2\nu_1=\nu_1\mu_2=0$,

\ \ \ \ $(\mu_1\nu_1)^m=(\mu_2\nu_2)^n, (\nu_1\mu_1)^m=(\nu_2\mu_2)^n, m,n\geqslant 1$.

$I_3:$ $\mu_1\nu_2=\nu_1\mu_1=\mu_2\nu_1=\nu_2\mu_2=0$,

\ \ \ \ $(\mu_1\nu_1\mu_2\nu_2)^n=(\mu_2\nu_2\mu_1\nu_1)^n, (\nu_1\mu_2\nu_2\mu_1)^n=(\nu_2\mu_1\nu_1\mu_2)^n, n\geqslant 1$.

\rule[0.5\baselineskip]{11.5cm}{1pt}

\ \ \ \ \ $Q: \xymatrix{\bullet \ar@<0.5ex>[r]^{\mu}\ar@(ld,lu)^{\alpha}&\bullet \ar@<0.5ex>[l]^{\nu}}$
\vspace{0.3cm}

$I_4:$ $\alpha\mu=\nu\alpha=0$, $\alpha^m=(\mu\nu)^n, m\geqslant 2, n\geqslant 1$.

$I_5:$ $\alpha^2=\nu\mu=0$, $(\alpha\mu\nu)^n=(\mu\nu\alpha)^n, n\geqslant 1$.

\rule[0.5\baselineskip]{8.5cm}{1pt}

\ \ \ \ \ $Q: \xymatrix{\bullet \ar@<0.5ex>[r]^{\mu}\ar@(ld,lu)^{\alpha}&\bullet \ar@<0.5ex>[l]^{\nu}\ar@(ru,rd)^{\beta}}$
\vspace{0.3cm}

$I_6:$ $\alpha\mu=\mu\beta=\beta\nu=\nu\alpha=0$, $\alpha^m=(\mu\nu)^n, \beta^r=(\nu\mu)^n,  m,r\geqslant 2, n\geqslant 1$.

$I_7:$ $\alpha^2=\nu\mu=\mu\beta=\beta\nu=0, (\alpha\mu\nu)^n=(\mu\nu\alpha)^n, \beta^m=(\nu\alpha\mu)^n, m\geqslant 2, n\geqslant 1$.

$I_8:$ $\alpha^2=\beta^2=\mu\nu=\nu\mu=0, (\nu\alpha\mu\beta)^n=(\beta\nu\alpha\mu)^n, (\alpha\mu\beta\nu)^n=(\mu\beta\nu\alpha)^n, n\geqslant 1$.

\rule[0.5\baselineskip]{14cm}{1pt}

\noindent In the above, we assume that $m,n,r\in \mathbb{N}$.
\end{proposition}

Then, we have the following observation.
\begin{proposition}\label{application}
Let $B_i$ be a two-point symmetric special biserial algebra. Then, $B_i$ is $\tau$-tilting finite if $i=1,4,5,6,7$; $\tau$-tilting infinite if $i=2,3,8$. Moreover, we have
\begin{center}
\renewcommand\arraystretch{1.1}
\begin{tabular}{|c|c|c|c|c|c|}
\hline
$B_i$&$B_1$&$B_4$&$B_5$&$B_6$&$B_7$\\ \hline
$\#\mathsf{s\tau\text{-}tilt}\ B_i$&\multicolumn{2}{c|}{6}&8&6&8\\ \hline
$Type$&\multicolumn{2}{c|}{$\mathcal{H}_{2,2}$}&$\mathcal{H}_{3,3}$&$\mathcal{H}_{2,2}$&$\mathcal{H}_{3,3}$\\
\hline
\end{tabular}.
\end{center}
\end{proposition}
\begin{proof}
One can easily check that $B_2$ and $B_3$ have $\Lambda_2$ as a quotient algebra, and $B_8$ has $\Lambda_5$ as a quotient algebra. Therefore, $B_2$, $B_3$ and $B_8$ are $\tau$-tilting infinite.

Then, we show the remaining case by case.

Case ($B_1$). If $n=1$, then $\mu\nu, \nu\mu \in C(B_1)$. If $n\geqslant 2$, then $\mu\nu+\nu\mu \in C(B_1)$. Both of them satisfy $\mathsf{s\tau\text{-}tilt}\ B_1\simeq \mathsf{s\tau\text{-}tilt}\ \left (B_1/<\mu\nu,\nu\mu>\right )\simeq \mathsf{s\tau\text{-}tilt}\ \Lambda_6$.

Case ($B_4$). If $n=1$, then $\alpha, \nu\mu \in C(B_4)$. If $n\geqslant 2$, then $\alpha, \mu\nu+\nu\mu \in C(B_4)$. Both of them satisfy $\mathsf{s\tau\text{-}tilt}\ B_4\simeq \mathsf{s\tau\text{-}tilt}\ \left (B_4/<\alpha,\mu\nu,\nu\mu>\right )\simeq \mathsf{s\tau\text{-}tilt}\ \Lambda_6$.

Case ($B_5$). If $n=1$, then $\mu\nu, \nu\alpha\mu\in C(B_5)$. If $n\geqslant 2$, then $\alpha\mu\nu+\mu\nu\alpha+\nu\alpha\mu \in C(B_5)$ and $\alpha\mu\nu\in C(\widetilde{B}_{5})$ such that $\mu\nu\in C (\widetilde{B}_{5}/<\alpha\mu\nu>)$, where
\begin{center}
$\widetilde{B}_{5}:= B_5/<\nu\alpha\mu,\alpha\mu\nu+\mu\nu\alpha>$.
\end{center}
Hence, $\mathsf{s\tau\text{-}tilt}\ B_5\simeq \mathsf{s\tau\text{-}tilt}\ \left (B_5/<\mu\nu,\nu\alpha\mu>\right )\simeq \mathsf{s\tau\text{-}tilt}\ \Lambda_8$.

Case ($B_6$). If $n=1$, then $\alpha, \beta\in C(B_6)$. If $n\geqslant 2$, then $\alpha, \beta, \mu\nu+\nu\mu \in C(B_6)$. Both of them satisfy $\mathsf{s\tau\text{-}tilt}\ B_6\simeq \mathsf{s\tau\text{-}tilt}\ \left (B_6/<\alpha,\beta,\mu\nu,\nu\mu>\right )\simeq \mathsf{s\tau\text{-}tilt}\ \Lambda_6$.

Case ($B_7$). If $n=1$, then $\beta,\mu\nu\in C(B_7)$. If $n\geqslant 2$, then $\beta, \alpha\mu\nu+\mu\nu\alpha+\nu\alpha\mu \in C(B_7)$ and $\alpha\mu\nu\in C(\widetilde{B}_{7})$ such that $\mu\nu\in C (\widetilde{B}_{7}/<\alpha\mu\nu>)$, where
\begin{center}
$\widetilde{B}_{7}:= B_7/<\beta,\nu\alpha\mu,\alpha\mu\nu+\mu\nu\alpha>$.
\end{center}
Thus, $\mathsf{s\tau\text{-}tilt}\ B_7\simeq \mathsf{s\tau\text{-}tilt}\ \left (B_7/<\beta,\mu\nu,\nu\alpha\mu>\right )\simeq \mathsf{s\tau\text{-}tilt}\ \Lambda_8$.
\end{proof}

Lastly, let $K$ be an algebraically closed field of characteristic not equal to 2. We consider the tame block algebras of Hecke algebras of classical type over $K$, which are classified in \cite{Ariki-tame-block}.

\begin{proposition}{\rm{(\cite[Theorem 2]{Ariki-tame-block})}}
Let $\Lambda$ be a tame block of Hecke algebras of classical type over $K$. Then, it is Morita equivalent to one of $D_i=KQ/I_i$ below.

\vspace{0.2cm}
$Q: \xymatrix{\bullet \ar@(ld,lu)^{\alpha}\ar@(ru,rd)^{\beta}}$\ \ \ \ \ \ \ \ \  $I_1:$ $\alpha^2=\beta^2=0, \alpha\beta=\beta\alpha$.

\rule[0.5\baselineskip]{9cm}{1pt}

$Q: \xymatrix{\bullet \ar@<0.5ex>[r]^{\mu}\ar@(ld,lu)^{\alpha}&\bullet \ar@<0.5ex>[l]^{\nu}}$\ \ \ \ \ \  $I_2:$ $\alpha\mu=\nu\alpha=0, \alpha^2=(\mu\nu)^2$.

\rule[0.5\baselineskip]{9cm}{1pt}

$Q: \xymatrix{\bullet \ar@<0.5ex>[r]^{\mu}\ar@(ld,lu)^{\alpha}&\bullet \ar@<0.5ex>[l]^{\nu}\ar@(ru,rd)^{\beta}}$

$I_3:$ $\alpha\mu=\mu\beta=\beta\nu=\nu\alpha=0, \alpha^2=\mu\nu, \beta^2=\nu\mu$.

$I_4:$ $\alpha\mu=\mu\beta=\beta\nu=\nu\alpha=0, \alpha^2=(\mu\nu)^2, \beta^2=(\nu\mu)^2$.

\rule[0.5\baselineskip]{9.7cm}{1pt}
\vspace{-0.2cm}
\end{proposition}

\begin{proposition}
Let $\Lambda$ be a tame block algebra of Hecke algebras of classical type over $K$. Then, $\Lambda$ is $\tau$-tilting finite. Moreover, we have
\begin{center}
\renewcommand\arraystretch{1}
\begin{tabular}{|c|c|c|c|c|}
\hline
$D_i$&$D_1$&$D_2$&$D_3$&$D_4$\\ \hline
$\#\mathsf{s\tau\text{-}tilt}\ D_i$&2&\multicolumn{3}{|c|}{6}\\ \hline
$Type$&$\mathcal{H}_{0,0}$&\multicolumn{3}{|c|}{$\mathcal{H}_{2,2}$}\\ \hline
\end{tabular}
\end{center}
\end{proposition}
\begin{proof}
Obviously, $D_1$ has only two basic support $\tau$-tilting modules, which are $D_1$ and 0. We can easily find the central elements for $D_2$, $D_3$ and $D_4$. Then, we can find that
\begin{center}
$\mathsf{s\tau\text{-}tilt}\ D_2\simeq \mathsf{s\tau\text{-}tilt}\ D_3\simeq \mathsf{s\tau\text{-}tilt}\ D_4\simeq \mathsf{s\tau\text{-}tilt}\ \Lambda_6$.
\end{center}
by applying Lemma \ref{centers}.
\end{proof}

\vspace{0.5cm}
\appendix
\section{Table T and Table W introduced in \cite{Han-wild-two-point}}
\begin{center}
\textbf{Table T\ \ \ \ \ \ \ \ \ }
\end{center}

$\xymatrix@C=1.5cm{\bullet\ar@<1.7ex>[r]^{\mu_1,\mu_2}\ar@<0.5ex>[r]&\bullet\ar@<1.7ex>[l]^{\nu_1,\nu_2}\ar@<0.5ex>[l]}$

\noindent(1) $\nu_1\mu_1=\nu_2\mu_2=(\ell_1\mu_1+\ell_2\mu_2)(k_1\nu_1+k_2\nu_2)=(\ell_3\mu_1+\ell_4\mu_2)(k_3\nu_1+k_4\nu_2)=0$,

where $k_1, k_2, k_3, k_4, \ell_1, \ell_2, \ell_3, \ell_4 \in K$ and $k_1k_4\neq k_2k_3, \ell_1\ell_4\neq \ell_2\ell_3$.

\noindent\rule[0.5\baselineskip]{16cm}{1pt}
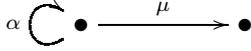
\begin{multicols}{2}
$\xymatrix@C=1.5cm{\bullet \ar[r]^{\mu}\ar@(dl,ul)^{\alpha }& \bullet }$
\vspace{0.3cm}

\noindent(2) $\alpha^6=\alpha^2\mu=0$;

\noindent\rule[0.5\baselineskip]{7cm}{1pt}

$
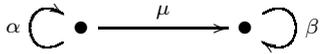
\xymatrix@C=1.5cm{\bullet \ar[r]^{\mu}\ar@(dl,ul)^{\alpha}& \bullet \ar@(ur,dr)^{\beta}}$\\

\noindent(3) $\alpha^2=\beta^2=0$;

\noindent(4) $\alpha^2=\beta^n=\mu\beta=0$, $n\geqslant 2, n\in \mathbb{N}$;

\noindent(5) $\alpha^m=\beta^n=\alpha\mu=\mu\beta=0$, $m, n\geqslant 2, m, n\in \mathbb{N}$;

\noindent(6) $\alpha^2=\beta^3=0, \alpha\mu=\mu\beta^2$;

\noindent(7) $\alpha^3=\beta^6=0, \alpha\mu=\mu\beta$;

\noindent(8) $\alpha^4=\beta^4=0, \alpha\mu=\mu\beta$;

\noindent\rule[0.5\baselineskip]{7cm}{1pt}

$
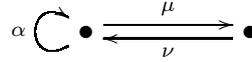
\xymatrix@C=1.5cm{\bullet \ar@<0.5ex>[r]^{\mu}\ar@(dl,ul)^{\alpha}& \bullet \ar@<0.5ex>[l]^{\nu}}$\\

\noindent(9) $\alpha^2=\mu\nu\mu=\nu\mu\nu=(\nu\alpha\mu)^n=0$, $n\geqslant 1, n\in \mathbb{N}$;

\noindent(10) $\alpha^3=\mu\nu=\nu\mu=\nu\alpha\mu=0$;

\noindent(11) $\alpha^3=\mu\nu, \nu\mu=\nu\alpha^2=\alpha^2\mu=0$;

\noindent(12) $\alpha^4=\mu\nu, \nu\alpha=\alpha^2\mu=0$;

\noindent(13) $\alpha^m=\nu\alpha=\alpha\mu=(\mu\nu)^n=0$, $m\geqslant 2, n\geqslant 1, m, n\in \mathbb{N}$;

\noindent(14) $\alpha^2=\mu\nu, \nu\alpha\mu=0$;

\noindent(15) $\alpha^3=\mu\nu, \nu\alpha=\alpha^2\mu=0$;

\noindent(16) $\alpha^3=\mu\nu, \nu\alpha=\nu\mu=0$;

\noindent\rule[0.5\baselineskip]{8cm}{1pt}
\end{multicols}
$
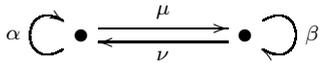
\xymatrix@C=1.5cm{\bullet \ar@<0.5ex>[r]^{\mu}\ar@(dl,ul)^{\alpha}&\bullet \ar@<0.5ex>[l]^{\nu}\ar@(ur,dr)^{\beta}}$
\vspace{0.3cm}

\noindent(17) $\alpha^2=\beta^2=\nu\mu=\mu\nu=0$;

\noindent(18) $\alpha^2=\beta^m=\nu\mu=\mu\beta=\beta\nu=(\nu\alpha\mu)^n=0$, $m\geqslant 2, n\geqslant 1, m, n\in \mathbb{N}$;

\noindent(19) $\alpha^m=\beta^n=(\nu\mu)^r=\alpha\mu=\nu\alpha=\mu\beta=\beta\nu=0$, $m,n\geqslant 2, r\geqslant 1, m, n, r\in \mathbb{N}$;

\noindent(20) $\alpha^2=\mu\nu, \beta^2=\nu\mu, \beta\nu=0, \alpha\mu=k\mu\beta$, $k\in K/\{0\}$;

\noindent(21) $\alpha^m=\beta^n=0, \beta^2=\nu\mu, \nu\alpha=\beta\nu$, $k_1\alpha^2=\mu\nu, \alpha\mu=k_2\mu\beta$,

\ \ \ $k_1,k_2\in K/\{0\}, m,n\geqslant 2, m, n\in \mathbb{N}$.

\noindent\rule[0.5\baselineskip]{16cm}{1pt}

\begin{center}
\textbf{Table W\ \ \ \ \ \ \ \ \ }
\end{center}
\begin{multicols}{2}
$\xymatrix@C=1.5cm{\bullet\ar@<1ex>[r]^{\mu_1,\mu_2,\mu_3}\ar[r]\ar@<-1ex>[r]&\bullet}$

\noindent(1) $KQ$;

\noindent\rule[0.5\baselineskip]{7cm}{1pt}

$\xymatrix@C=1.5cm{\bullet\ar@<1ex>[r]^{\mu_1,\mu_2}\ar[r]&\bullet\ar@<1ex>[l]^{\nu}}$

\noindent(2) $\mu_1\nu=\mu_2\nu=0$;

\noindent\rule[0.5\baselineskip]{8cm}{1pt}
\end{multicols}
\vspace{-0.3cm}

$\xymatrix@C=1.5cm{\bullet\ar@<1.7ex>[r]^{\mu_1,\mu_2}\ar@<0.5ex>[r]&\bullet\ar@<1.7ex>[l]^{\nu_1,\nu_2}\ar@<0.5ex>[l]}$

\noindent(3) $\nu_2\mu_1=\nu_1\mu_2, \mu_1\nu_1=\mu_2\nu_1=\mu_1\nu_2=\mu_2\nu_2=\nu_1\mu_1=0$;

\noindent\rule[0.5\baselineskip]{12cm}{1pt}
\begin{multicols}{2}
$\xymatrix@C=1.5cm{\bullet\ar[r]^{\mu}\ar@(ul,ur)^{\alpha_1}\ar@(dl,dr)_{\alpha_2}& \bullet }$

\noindent(4) $\alpha_1^2=\alpha_2^2=\alpha_1\alpha_2=\alpha_2\alpha_1=\alpha_1\mu=\alpha_2\mu=0$;

\noindent\rule[0.5\baselineskip]{7cm}{1pt}

$\xymatrix@C=1.5cm{\bullet\ar@(ul,ur)^{\alpha}\ar@<0.7ex>[r]^{\mu_1}\ar@<-0.7ex>[r]_{\mu_2}&\bullet}$\\

\

\noindent(5) $\alpha^2=\alpha\mu_1=\alpha\mu_2=0$;

\noindent\rule[0.5\baselineskip]{8cm}{1pt}
\end{multicols}
\begin{multicols}{2}
$\xymatrix@C=1.5cm{\bullet \ar[r]^{\mu}\ar@(dl,ul)^{\alpha }& \bullet }$
\vspace{0.5cm}

\noindent(6) $\alpha^7=\alpha^2\mu=0$;

\noindent(7) $\alpha^4=\alpha^3\mu=0$;

\noindent\rule[0.5\baselineskip]{7cm}{1pt}

$\xymatrix@C=1.5cm{\bullet \ar[r]^{\mu}\ar@(dl,ul)^{\alpha}& \bullet \ar@(ur, dr)^{\beta}}$\\

\noindent(8) $\alpha^2=\beta^3=\alpha\mu=0$;

\noindent(9) $\alpha^3=\beta^3=\alpha\mu=\mu\beta^2=0$;

\noindent(10) $\alpha^2=\beta^4=\alpha\mu=\mu\beta^2=0$;

\noindent(11) $\alpha^2=\beta^3=\alpha\mu\beta=\mu\beta^2=0$;

\noindent(12) $\alpha^4=\beta^5=\mu\beta^2=0, \alpha\mu=\mu\beta$;

\noindent(13) $\alpha^3=\beta^7=\mu\beta^2=0, \alpha\mu=\mu\beta$;

\noindent\rule[0.5\baselineskip]{7cm}{1pt}

$\xymatrix@C=1.5cm{\bullet \ar@<0.7ex>[r]^{\mu}\ar@(dl,ul)^{\alpha}& \bullet \ar@<0.7ex>[l]^{\nu}}$
\vspace{0.3cm}

\noindent(14) $\alpha^3=\mu\nu=\nu\mu=\alpha^2\mu=0$;

\noindent(15) $\alpha^3=\mu\nu=\alpha\mu=0$;

\noindent(16) $\alpha^3=\mu\nu=\nu\alpha\mu=\nu\alpha^2=\alpha^2\mu=0$;

\noindent(17) $\alpha^4=\mu\nu=\nu\mu=\alpha\mu=\nu\alpha^3=0$;

\noindent(18) $\alpha^4=\mu\nu=\nu\mu=\nu\alpha\mu=\nu\alpha^2=\alpha^2\mu=0$;

\noindent(19) $\alpha^5=\mu\nu=\nu\mu=\nu\alpha=\alpha^2\mu=0$;

\noindent(20) $\alpha^2=\nu\alpha=\nu\mu\nu=\alpha\mu\nu=0$;

\noindent(21) $\alpha^2=\nu\alpha=\mu\nu\mu=0$;

\noindent(22) $\alpha^3=\nu\mu=\nu\alpha=\alpha\mu\nu=\alpha^2\mu=0$;

\noindent(23) $\alpha^2=\mu\nu, \alpha^3=\alpha^2\mu=0$;

\noindent(24) $\alpha^3=\mu\nu, \alpha^4=\nu\mu=\nu\alpha\mu=\nu\alpha^2=0$;

\noindent\rule[0.5\baselineskip]{8cm}{1pt}
\end{multicols}

$\xymatrix@C=1.5cm{\bullet \ar@<0.7ex>[r]^{\mu}\ar@(dl,ul)^{\alpha}&\bullet \ar@<0.7ex>[l]^{\nu}\ar@(ur,dr)^{\beta}}$
\vspace{0.3cm}

\noindent(25) $\alpha^3=\beta^2=\nu\mu=\mu\nu=\nu\alpha=\mu\beta=\beta\nu=\alpha^2\mu=0$;

\noindent(26) $\alpha^2=\beta^2=\nu\mu=\alpha\mu=\nu\alpha=\beta\nu=0$;

\noindent(27) $\alpha^2=\mu\nu, \beta^2=\nu\mu=\alpha\mu=\mu\beta=\beta\nu\alpha=0$;

\noindent(28) $\alpha^2=\mu\nu, \beta^2=\nu\mu=\alpha\mu=\beta\nu=0$;

\noindent(29) $\alpha^2=\mu\nu, \beta^2=\nu\mu=\nu\alpha=\mu\beta=0$;

\noindent(30) $\alpha^2=\mu\nu, \beta^2=\nu\mu=\nu\alpha=\beta\nu=\alpha\mu\beta=0$;

\noindent(31) $\alpha\mu=\mu\beta, \alpha^2=\beta^3=\mu\nu=\nu\alpha=\beta\nu=\mu\beta^2=0$;

\noindent(32) $\alpha\mu=\mu\beta, \alpha^2=\beta^2=\nu\alpha=\beta\nu=\mu\nu\mu=\nu\mu\nu=0$;

\noindent(33) $\alpha\mu=\mu\beta, \alpha^3=\beta^3=\nu\mu=\mu\nu=\nu\alpha=\beta\nu=\mu\beta^2=\alpha^2\mu=0$;

\noindent(34) $\alpha\mu=\mu\beta, \alpha^3=\beta^2=\nu\mu=\nu\alpha=\beta\nu=\alpha^2\mu=0$;

\noindent\rule[0.5\baselineskip]{16cm}{1pt}

\ \\

Department of Pure and Applied Mathematics, Graduate School of Information Science and Technology, Osaka University, 1-5 Yamadaoka, Suita, Osaka, 565-0871, Japan.

\emph{Email address}: \texttt{q.wang@ist.osaka-u.ac.jp}
\end{spacing}
\end{document}